\documentclass[a4paper,12pt]{report}

\newcommand{\usecustompackage}[2][]{%
	\IfFileExists{Packages/#2.sty}{%
		\usepackage[#1]{./Packages/#2}%
	}{%
		\usepackage[#1]{#2}%
	}%
}

\usepackage[T1]{fontenc}
\usepackage[utf8]{inputenc}
\usepackage[english]{babel}
\usepackage[left=2.2cm, right=2.2cm, top=2.2cm, bottom=2.2cm]{geometry}

\usecustompackage{common_packages}

\usepackage{caption}

\usecustompackage{macros_generic}
\usecustompackage[cnt1=chapter, cnt2=section]{macros_theorems}
\usecustompackage{macros_tikz}

\fancypagestyle{nonumber}{
	\fancyhf{}
}

\fancypagestyle{numbered}{
	\fancyhf{}
	\fancyhead[L]{}
	\fancyhead[C]{}
	\fancyhead[R]{}
	\fancyfoot[L]{}
	\fancyfoot[C]{\arabic{page}}
	\fancyfoot[R]{}
}

\newproof{pfsketch}{A bizonyítás vázlata}{Sketch of the proof}{\bfseries}{}

\newcommand{\nt}[1]{\It{#1}}

\begin{document}

\pagestyle{nonumber}
\begin{center}
	\vspace*{0.8cm}
	{\large\textbf{MASTER'S THESIS}}\\
	\vspace{3cm}

	\Bf{\large Probabilistic formulation of the Hadwiger--Nelson problem}\\
	 \vspace{4cm}

	\Bf{\large Péter Ágoston}\\
	\vspace{0.2cm}
	\Bf{MSc in Mathematics}\\
	\vspace{2cm}

	\begin{tabular}{rl}
		\Bf{\large Supervisor:} & Dömötör Pálvölgyi\\
			& assistant professor\\
			& ELTE TTK Institute of Mathematics\\
			& Department of Computer Science\\
	\end{tabular}
	\vspace{1cm}
	\begin{figure}[h!]
		\centering
		\includegraphics[width=50mm]{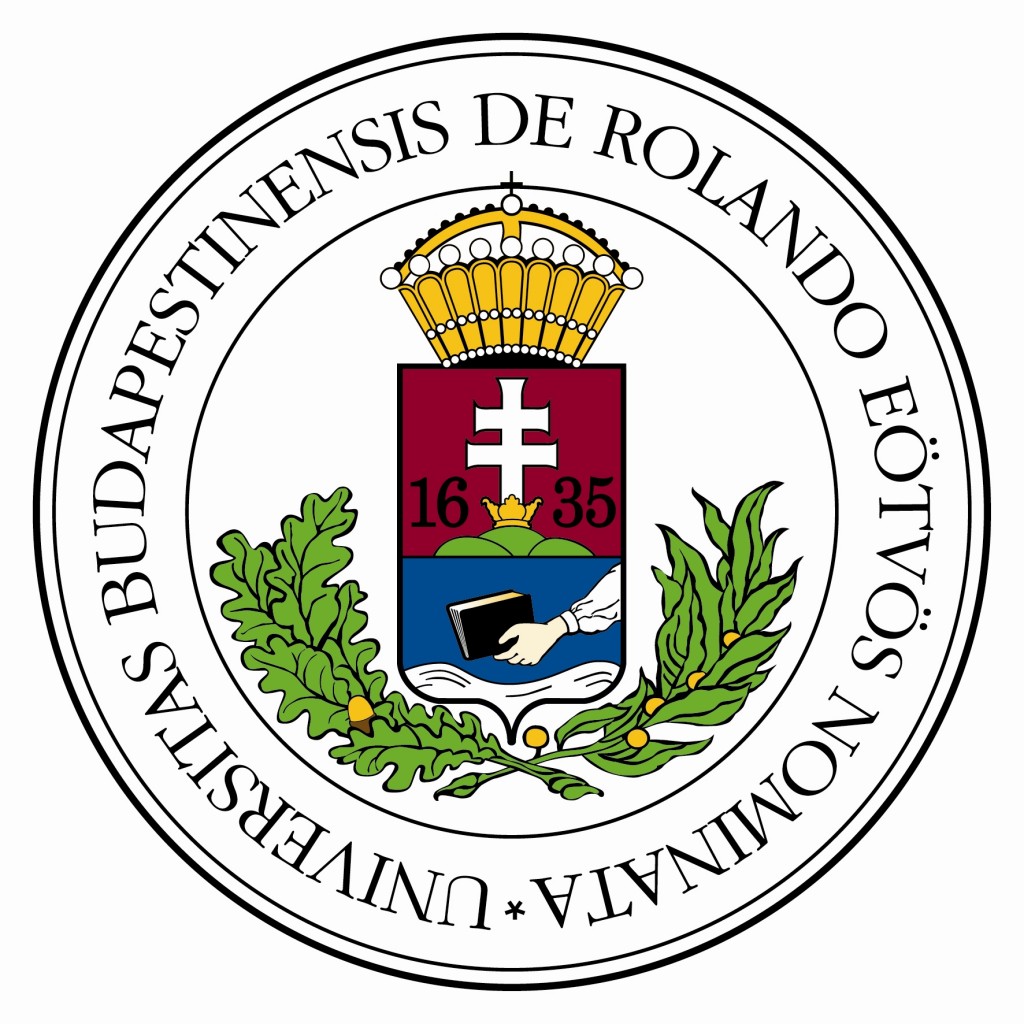}
	\end{figure}

	\Bf{\large ELTE}\\
	\Bf{\large 2019}\\
\end{center}
\pagebreak

\pagestyle{numbered}
\tableofcontents

\chapter*{Introduction}
\addcontentsline{toc}{chapter}{Introduction}

In 1950, Edward Nelson posed the question of determining the chromatic number of the graph on the plane formed by its points as vertices and having edges between the pairs of points with distance $1$, which is now simply known as the chromatic number of the plane. It soon became obvious that at least $4$ colours are needed, which can be easily seen thanks to a graph with $7$ vertices, the so called Moser spindle, found by William and Leo Moser. There is also a relatively simple colouring by John R. Isbell, which shows that the chromatic number is at most $7$. And despite the simplicity of the bounds, these remained the strongest ones until $2018$, when biologist Aubrey de Grey found a subgraph which was not colourable with $4$ colours. This meant that the chromatic number of the plane is at least $5$. The proof used a computer program to verify that there is no $4$-colouring of the graph, so it was a natural question to find a smaller subgraph which still cannot be coloured with $4$ colours. A Polymath project was launched which had the goal to simplify this bounding graph. One of the methods used is a probabilistic one, which, given a colouring, assigns probabilities to every distance, so that it tells the probability of two points with the given distance being monochromatic and we can make bounds for these probabilities. The results achieved with this method may also be helpful in related problems, for example questions about graphs which also include pairs of points with some other distances.

The first part of my thesis describes basic properties and the possible edge count of unit distance graphs, graphs which can be realized in the plane so that the distance of their adjacent vertices is $1$, as the Hadwiger--Nelson problem is based on these graphs. As a new result, this chapter also contains an improvement with a constant factor of the currently known best bound of the maximal edge count of a unit distance graph in Theorem \thmref{unfelsobecsles}

In Chapter $2$, we introduce $(1,d)$-graphs, a similar concept to unit distance graphs, but here we allow two distances to occur between the vertices. These graphs can be used in the probabilistic approach of the Hadwiger--Nelson problem.

Chapter $3$ deals with the chromatic number of the plane, so with the Hadwiger--Nelson problem in general.

Chapter $4$ contains estimates from the probabilistic approach for the Hadwiger--Nelson problem, including several results obtained in joint work with my supervisor.

Finally in Chapter $5$ gives an outlook to related problems on spheres and in other dimensions.

This is a corrected and slightly extended version of my thesis, the original version can be found here: \href{https://web.cs.elte.hu/blobs/diplomamunkak/msc\_mat/2019/agoston\_peter.pdf}{https://web.cs.elte.hu/blobs/diplomamunkak/msc\_mat/2019/agoston\_peter.pdf}

\chapter*{Acknowledgements}

I would like to thank my supervisor Professor Dömötör Pálvölgyi for introducing me to the topic of my thesis. In particular I am thankful for his patience and persistence with which he helped me enter deeper into the field through many inspiring hours of joint work.

Special thanks go to Professor Heiko Harborth from TU Braunschweig who sent me an otherwise unattainable work, the thesis of his former student Carsten Schade.

Finally, I express my gratitude to my family for their constant support and advice during the preparation of the thesis.

\chapter{Unit distance graphs}

\section{Introduction}

First we list some elementary properties of our main objects of investigation, the unit distance graphs.

\begin{defn}
	We call a graph a \nt{unit distance graph} (UDG) if its vertices can be represented by distinct points in
	the plane so that two points have distance $1$ if they are connected with an edge. Such a representation
	is called a unit distance representation (UDR) of the graph.
\end{defn}

\begin{defn}
	We call a graph a \nt{faithful unit distance graph} (FUDG) if its vertices can be represented by distinct points
	in the plane so that two points have distance $1$ if and only if they are connected with an edge. Such a
	representation is called a faithful unit distance representation (FUDR) of the graph.
\end{defn}

\begin{ex}[label=notfaithful]
	Although the UDG in Figure \ref{udgnemfudg} is a unit distance graph, it is not faithful: the triangles:
	$ABC$, $ACD$, $ADE$, $AEF$ are equilateral in every UDR and they only can be placed in one way up
	to isometry, so the position of these $6$ points is unique up to isometry and we get that although $AG$
	is not an edge, the only place $G$ can be is at distance $1$ from $A$, which means that there is no faithful
	unit distance representation of this graph. However if we add edge $AG$ as in Figure \ref{fudgreszenemfudg},
	we get a FUDG since all the non-edges have length $\sqrt{3}$ or $2$.
\end{ex}

\begin{center}
	\begin{minipage}{.49\textwidth}
		\centering
		\input{Figures/UDGaminemFUDG}
		\captionof{figure}{}
		\label{udgnemfudg}
	\end{minipage}
	\begin{minipage}{.49\textwidth}
		\centering
		\input{Figures/FUDGaminekegyreszenemFUDG}
		\captionof{figure}{}
		\label{fudgreszenemfudg}
	\end{minipage}%
\end{center}

\begin{prop}[label=udgsubgraphfudg]
	For every unit distance graph $G$, there is a faithful unit distance graph, which has $G$ as a spanning subgraph.
\end{prop}

\begin{pf}
	Since $G$ is a UDG, we have a UDR for it, so if we take the $1$ distances and draw edges between all of the
	corresponding pairs of vertices, the UDR for $G$ will be a FUDR for this new graph.
\end{pf}

\begin{defn}
	A property of a graph is \nt{hereditary} if it is also true for all of its induced subgraphs.
\end{defn}

\begin{defn}
	A property of a graph is \nt{monotone} if it is also true for all of its subgraphs.
\end{defn}

It is clear that if a property is monotone, it is also hereditary.

\begin{prop}[label=monotone]
	The unit distance graph property is monotone, so it is also hereditary.
\end{prop}

\begin{pf}
	If there is a unit distance representation of a graph $G$ in the plane, its respective subset is also a unit
	distance representation for an arbitrary subgraph: the points in the subset are still distinct and the edges
	were all edges in $G$, so the points belonging to its end vertices have distance $1$.
\end{pf}

\begin{prop}
	The faithful unit distance graph property is hereditary, but not monotone.
\end{prop}

\begin{pf}
	If we have a FUDG, then the respective subset of its FUDR will be a FUDR for an arbitrary subgraph,
	since all the $1$ distances will be edges and all the other distances will be non-edges, so the FUDG
	property is hereditary. Also, Example \thmref{notfaithful} shows that this property is not monotone.
\end{pf}

\section{Upper bound on the number of edges}

\begin{defn}
	Let us denote the maximal number of edges in a unit distance graph with $n$ vertices by $u(n)$.
\end{defn}

\begin{prop}[label=maxudgfudg]
The maximal possible number of edges in a faithful unit distance graph with $n$ vertices is also $u(n)$.
\end{prop}

\begin{pf}
All faithful unit distance graphs are unit distance graphs too, so there cannot be a larger faithful unit distance graph with $n$ vertices than the maximal unit distance graph(s).

Also if a maximal unit distance graph with $n$ vertices would not be a faithful unit distance graph that would mean that it is the subgraph of a larger unit distance graph with $n$ vertices because of Proposition \thmref{udgsubgraphfudg}, which is a contradiction, so all of the maximal unit distance graphs with $n$ vertices are also faithful unit distance graphs, which means that we are done.
\end{pf}

\begin{prop}
$u(2n)\ge2\cdot u(n)+n$ for $n\ge0$.
\end{prop}

\begin{pf}
We can construct a unit distance graph with $2n$ vertices and $2\cdot u(n)+n$ edges: just take two copies of the same UDR of the same maximal unit distance graph with $n$ vertices and place them so that they are the translates of each other with a vector of length $1$ that does not coincide with any of the vectors defined by points in one of the copies. That makes $n+n=2n$ vertices, $u(n)$ edges in the first copy, $u(n)$ edges in the second copy and at least $n$ edges between them.
\end{pf}

\begin{prop}[label=minkowski]
$u(ab)\ge a\cdot u(b)+b\cdot u(a)$.
\end{prop}

\begin{pf}
The proof is similar to the one above: if we have a UDG with $a$ vertices and $u(a)$ edges and another with $b$ vertices and $u(b)$ edges, we can take some UDR ($P_1$ and $P_2$) for both of them and rotate one of them so that we cannot choose a vector between two points of $P_1$ and another between two points of $P_2$ so that their sum is $0$. We can do that since for any two vectors, there is at most one rotation, which makes their sum $0$ and there is a finite number of possible pairs of vectors. Now, if we take the Minkowski sum of the two UDRs, that will be a point set with $a\cdot u(b)+b\cdot u(a)$ pairs of points with distance $1$.
\end{pf}

\begin{prop}[label=udgedgedensity]
The maximal possible edge density of a unit distance graph with $n$ vertices is monotonically decreasing, or more specifically, $u(n)\le{n\over n-2}\cdot u(n-1)$ for $n\ge 1$.
\end{prop}

\begin{pf}
In Proposition \thmref{monotone} we have seen that the unit distance graph property is monotone. From here it is enough to prove the following lemma:

\begin{lem}[label=monedgedensity]
For a monotone property on graphs, the maximal possible edge density of a graph with $n$ vertices satisfying the property is monotonically decreasing, or equivalently for the maximum number $m(n)$ of the edges, $m(n)\le{n\over n-2}\cdot m(n-1)$ for $n\ge 1$.
\end{lem}

\begin{pf}
Let us take a maximal graph $G$ with $n$ vertices ($n\ge1$) which has the given property and let us take all of its subgraphs with $n-1$ vertices. There are $n$ such subgraphs and if we add up the number of their edges, all of the edges of $G$ have been counted $n-2$ times (an edge is missed exactly when one of its endpoints is left out) and since all of the subgraphs have the given property, the sum is at most $n\cdot m(n-1)$, which implies $(n-2)\cdot m(n)\le n\cdot m(n-1)$. So for the maximal edge densities:
$$d(n)={m(n)\over{n\choose2}}\le{n\over n-2}\cdot{m(n-1)\over{n\choose2}}={m(n-1)\over{n-1\choose2}}=d(n-1)$$
\end{pf}
\end{pf}

\begin{defn}
Let us draw a (not necessarily simple) graph in the plane so that the images of the edges do not go through the images of the vertices (except for the endpoints of that certain edge, but even then, they only end there). The \nt{crossing number} of a graph is the minimum number of points where the images of the edges intersect (counted with multiplicity, so for example if $n$ edges intersect in the same point, that counts as ${n\choose2}$ crossings). (Note: alternatively, we can define it so that we do not allow more than $2$ edges to intersect in one point, since such a planar embedding can be achieved by a very small change to the images of the edges.) We denote the crossing number of a graph $G$ with $cr(G)$.
\end{defn}

\begin{defn}
Let us draw two (not necessarily simple) graphs in the plane on the same vertex set so that the images of the edges do not go through any vertices and the set of their edges is disjoint. The \nt{crossing number of the two graphs} is the minimum number of points (again, counted with multiplicity) where the images of the edges of the first graph intersect the images of the edges of the second one. Let us denote the crossing number of graph $G$ and $H$ with $cr(G,H)$.
\end{defn}

The so called crossing number inequality or crossing lemma has several versions, the next two statements are the two best currently known versions of it.

\begin{thm}[label=regicrossinglemma, description={Pach, Radoičić, Tardos, Tóth, 1997}]{\rm{\cite{prtt}}} If a simple graph has $n$ vertices and $e$ edges so that $e>4n$, then the crossing number of the graph is at least $\displaystyle{e^3\over64n^2}$.
\end{thm}

\begin{thm}[label=ujcrossinglemma, description={Ackermann, 2013}]{\rm{\cite{a}}} If a simple graph has $n$ vertices and $e$ edges so that $e>7n$, then the crossing number of the graph is at least $\displaystyle{e^3\over29n^2}$.
\end{thm}

\begin{thm}[label=unfelsobecsles]
$u(n)=O(n^{4/3})$. More specifically, we will prove that the maximum number of edges is at most $\sqrt[3]{{2\over3}{\left({2+\sqrt{3}}\right)}\cdot29}\cdot n^{4/3}$ for large enough $n$'s.
\end{thm}

The first proof for $O(n^{4/3})$ is by Spencer, Szemerédi and Trotter (1984) {\rm{\cite{sszt}}}. The proof below is an improvement of the one given by László Székely {\rm{\cite{sz1}}}, {\rm{\cite{sz2}}}, which achieves $8n^{4/3}$, but instead of using the older Theorem \thmref{regicrossinglemma} it uses the newer, improved version Theorem \thmref{ujcrossinglemma} and also has some other improvements.

\begin{pf}
Let us take a planar set $P$ of $n$ points, which is a UDR for the UDG with maximal edge count $q$ on $n$ vertices. Let us denote the corresponding maximal UDG (which is clearly a FUDG, as it was shown in the proof of Proposition \thmref{maxudgfudg}) by $G(P)$.

\begin{lem}[label=egyseggraflemma]
If $n\ge 3$, then all vertices of $G(P)$ have degree at least $2$, and if $n\ge 7$, then there is at least one unit distance graph with the same vertex and edge number as $P$, for which there is at most one vertex with degree $2$.
\end{lem}

\begin{pf}
If we have a maximal construction, there is an edge $AB$, which is not the side of two distinct equilateral triangles with vertices in $P$ (so out of the two possible triangles $ABC$ and $ABD$, at least one (say, $ABC$) has its third vertex ($C$) outside of $P$), otherwise we could find an infinite grid of equilateral triangles, since $G(P)$ is a FUDG, and every 1 distance also means an edge, where we can also place an equilateral triangle, whose third vertex also should be in $P$ by the assumption, etc. Let us choose a point $C$ so that it forms an equilateral triangle with side length $1$ with two vertices in $P$. If there is a vertex in $G(P)$, with degree less than $2$, then we can replace it by $C$, so its degree will be at least $2$ and the number of edges have increased, while the number of vertices does not change, so the original graph was not maximal. Thus the first part of the statement is true.

For the second part, let us suppose that for $n\ge 7$, we have a UDG with maximal edge count, which has more than $1$ edges with degree $2$. Now let us repeat the following algorithm: if we find a vertex with degree at most $2$, we delete it along with its neighbouring edges, and we do this until after $k$ deleted vertices there are no vertices left with degree at most $2$. Let us now suppose that no vertices remain after this procedure. This would mean that in the last step, we have deleted $0$ edges, in the second to last we have deleted at most $1$, and before that, in every step, we have deleted at most $2$ edges. So the original graph had at most $2n-3$ edges, which is impossible since it was a UDG with maximal edge count on $n$ vertices, and there exists a UDG with $n$ vertices and $2n-2$ edges for every $n\ge 7$, similar to the one in Figure \ref{2nminusz2el}:

\begin{center}
	\input{Figures/2nminusz2elugraf}
	\captionof{figure}{}
	\label{2nminusz2el}
\end{center}

So the previous assumption led us to a contradiction, hence the remaining point set $P'$ is not empty. So $G(P')$ has at least $1$ vertex and since all of its vertices have degree at least $3$, it has more than $0$ edges. So we may assume that it has an edge, which has an angle $+60^\circ$ with a directed horizontal line pointing to the right. Now let us make a triangular chain to the left as seen in Figure \ref{3szoglanc} (where the white points (and $C$) are in a triangular chain starting from the edge $AB$) by putting vertices according to the grid from the right to the left (and drawing all possible edges from them) except if a point is already in $P'$ (like $C$ in the Figure), when we simply skip it and move to the next point. (In the figure, the order of the procedure is represented by the numbering.) In every step, we draw at least $2$ edges, since all of the points have two neighbours to their right which are in the chain (if we include the starting edge). We do this procedure until we added $k$ points and get a graph $P''$. The procedure ends at some point, as $P'$ has finitely many points and $k$ is also finite. In the deleting algorithm, we deleted at most $2$ edges in every step, thus at most $2k$ edges were deleted, so $P''$ has at least as many edges as $G(P)$, since we have added at least $2k$ edges. And $P''$ has only one vertex with degree $2$: the endpoint of the chain. All the other points in the chain (except for the rightmost edge that is in $P'$) have degree at least $3$: they have at least $2$ neighbours to the right and at least $1$ to the left. The points in $P'$ already had degree more than $2$, so they still have. And since all points in $P''$ are part of either the triangle chain (without the rightmost edge) or $P'$, it means that we found a unit distance graph with the same number of vertices ($\abs{V(P)}=\abs{V(P')}+k=\abs{V(P'')}$), which has at most $1$ vertex with degree $2$ (it has exactly one such vertex if the original one had at least one).

\begin{center}
	\input{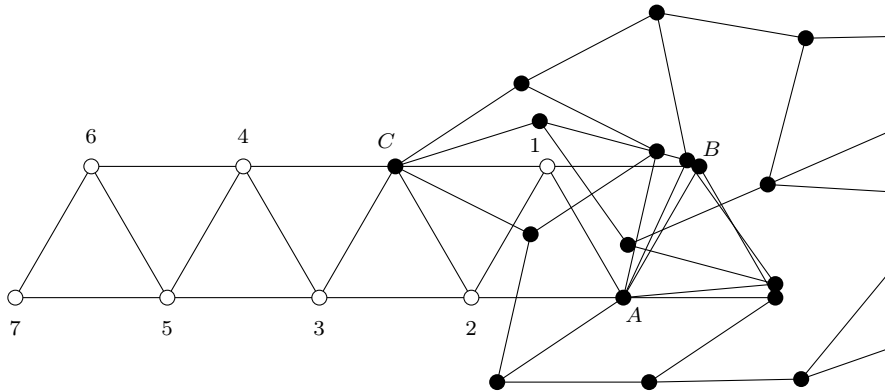}
	\captionof{figure}{An example graph: black points denote the vertices of $P'$ and the white ones are those added to it to create $P''$.}
	\label{3szoglanc}
\end{center}

(Note: after the end of the process, we can see that if there actually would be a point in the triangle chain, that was in $P'$, while the previous one was not, that is actually impossible as it contradicts to the assumption of $P$ being maximal, since then together with the previous point, we created at least $3$ edges, so there were more edges created in the creation process, as deleted in the deletion process.)
\end{pf}

From now on, we look only at the case, when $n\ge 7$ and we assume that $P$ was chosen in a way that it has no vertex with degree at most $1$ and has at most $1$ vertex with degree $2$ (according to the lemma, the first part is always true and the second part can be achieved by choosing $P''$ instead of the original $P$) and let $V$ be the potentially existing vertex with degree $2$.

Let us define a graph $G$ as follows:

Its vertices are the points of $P$ and its edges are defined by arcs on the unit circles around the points of $P$, whose endpoints are two neighbouring points on the circles (as we have seen above, all of the unit circles around the points of $P$ have at least two points from $P$ on them).

$\abs{E(G)}=2q$, since for every point of $P$, there are exactly as many arcs on the circle around it as the degree of the corresponding point in $G(P)$, so each edge is counted twice.

$G$ does not have loops, since we proved above that every vertex has at least $2$ neighbours, so every circle passes through at least $2$ points of $P$.

Multiple edges can occur, but there is at most one pair of points, which are connected on the same circle because of the assumption from the lemma (there is at most one vertex in $G(P)$, which has degree $2$ and for all the others, there cannot be two shortest arcs between two points). So if we leave out one of these two edges, we get that $E(G)\le2q-1$. The other possibility of having multiple edges is that there are two circles passing through both points, but there cannot be more, so know all edges occur with multiplicity at most $2$.

Also, the crossing number of $G$ is at most $2\cdot{n\choose2}=n\cdot(n-1)$, since all pairs of circles intersect each other in at most $2$ points, and that is the only way for two edges in $G$ to cross.

Let us now denote the number of edges with no duplicates by $e$ and the number of edges with one duplicate by $f$. Let us also delete the possibly existing edge with $3$ duplicates from $H$. Let us also call the graph induced by the $e$ unique edges as $H'$, the graph induced by the $f$ non-unique edges as $H''$. Let us also suppose that $e$ and $f$ are both larger than $7n$. This way, if we choose $H$ randomly, then the expected value of $cr(H')+2cr(H',H'')+4cr(H'')$ is at most $cr(G)\le n\cdot(n-1)=n^2-n$, so if we choose it appropriately, $cr(H')+2cr(H',H'')+4cr(H'')\le cr(G)$. Let us now call $E=\left(cr(H')+cr(H',H'')\right)\cdot29n^2$ and $F=cr(H'')\cdot29n^2$. This way we know that $E+4F\le 29\cdot(n^4-n^3)$, $e\le\sqrt[3]{E}$, $f\le\sqrt[3]{F}$, $e+f\le\sqrt[3]{E+F}$ and $2q\ge e+2f$ (the second, the third and the fourth inequality is due to Theorem \thmref{ujcrossinglemma} applied to $H'$, $H''$ and $H'+H''$). Let us now call $E+4F=D$. This means that $D\le cr(G)\cdot29n^2\le29\cdot(n^4-n^3)$. Then $e+2f=(e+f)+f\le\sqrt[3]{E+F}+\sqrt[3]{F}$ (and the bound is tight, since from the inequalities, it is possible that $e=\sqrt[3]{E+F}-\sqrt[3]{F}$, while $f=\sqrt[3]{F}$, since $\sqrt[3]{E+F}\le\sqrt[3]{E}+\sqrt[3]{F}$). But then if we write $x={E\over D}$, then $E=x\cdot D$ and $F={(E+4F)-E\over4}={D-xD\over4}={1-x\over4}\cdot D$. Now we have to calculate the maximum of $\sqrt[3]{E+F}+\sqrt[3]{F}=\sqrt[3]{x\cdot D+{1-x\over4}\cdot D}+\sqrt[3]{{1-x\over4}\cdot D}=\left(\sqrt[3]{x+{1-x\over4}}+\sqrt[3]{{1-x}\over4}\right)\cdot\sqrt[3]{D}$ for constant $D$ and $0\le x\le D$, so it is enough to maximalize $\sqrt[3]{x+{1-x\over4}}+\sqrt[3]{{1-x}\over4}$. And the maximum is $\sqrt[3]{{2\over3}\left(2+\sqrt{3}\right)}$ {\rm{\cite{wa1}}}. And from here the maximum of $e+2f$ is $\sqrt[3]{{2\over3}\left(2+\sqrt{3}\right)}\cdot\sqrt[3]{D}\le\sqrt[3]{{2\over3}\left(2+\sqrt{3}\right)}\cdot\sqrt[3]{29\cdot(n^4-n^3)}$

So $2q-1\le\sqrt[3]{{2\over3}\cdot\left(2+\sqrt{3}\right)\cdot29\cdot(n^4-n^3)}\Rightarrow q\le{\sqrt[3]{{2\over3}\cdot\left(2+\sqrt{3}\right)\cdot29}\over2}\cdot n^{4\over3}\approx2.082\cdot n^{4\over3}$.

Now let us solve the case, when $e\le 7n$ or $f\le 7n$. If $e\le 7n$, but $f>7n$ then $2q-1\le 7n+2\cdot\sqrt[3]{29\over4}\cdot n^{4\over3}$ with a similar crossing lemma argument if we choose $H''$ to have at most ${1\over4}$ of the crossings, and this is asymptotically smaller than $4\cdot n^{4\over3}$, so $q$ is asymptotically smaller than $2\cdot n^{4\over3}$. If only the second is true, then $2q-1\le 14n+\sqrt[3]{29}\cdot n^{4\over3}$ similarly, which is also asymptotically smaller than $2\cdot n^{4\over3}$, so $q$ is asymptotically smaller than $2\cdot n^{4\over3}$. If $e\le 7n$ and $f\le 7n$, then the edge count is trivially smaller asymptotically than $2\cdot\sqrt[3]{29n^4}$. So we are done with these cases.
\end{pf}

\begin{thm}[description={Erdős}]
{\rm{\cite{bmp}}} $u(n)>ne^{c\log{n}\over\log{\log{n}}}$ for some $c$.
\end{thm}

\begin{pfsketch}We can take appropriate scaled sections of the integer lattice, which give this bound:

For every number $r$ and a point $P$ of the integer lattice, the number of points with distance $\sqrt{r}$ from $P$ equals the number of integer pairs $a$ and $b$, whose squares sum up to $r$ (we can do a bijection between such pairs and points $P+(a,b)$). And for the product $r_k$ of the first $k$ primes, which are $\equiv1\pmod 4$, there are enough vectors with length $\sqrt{r_k}$ among the points of the integer lattice so that we can get a bound $u(n)\ge\left(ne^{c\log{n}\over\log{\log{n}}}\right)$.
\end{pfsketch}

\begin{con}[description={Erdős}]
{\rm{\cite{bmp}}} $u(n)\le O\left(n^{1+\varepsilon}\right)$ for all $\varepsilon>0$. A stronger version is that $u(n)\le O{\left(ne^{c\log{n}\over\log\log{n}}\right)}$.
\end{con}

\begin{prop}
[label=k4]{$K_4$ is not a unit distance graph.}
\end{prop}

\begin{pf}
If there would be a UDR for $K_4$, then all four points in it should be of distance $1$ from each other, which would mean that the first two are on a segment of length $1$ and the other two are the vertices of equilateral triangles which both have this segment as one of their sides, but then the only possible point set is the one in Figure \ref{4csucs5el}, where the distance of $C$ and $D$ is $\sqrt{3}$ and not $1$, so it is not a good point set.
\end{pf}

\begin{center}
	\input{Figures/maximalis4csucsuUDG}
	\captionof{figure}{}
	\label{4csucs5el}
\end{center}

\begin{prop}[label=k23]
$K_{2,3}$ is not a unit distance graph.
\end{prop}

\begin{pf}
If there would be a UDR for $K_{2,3}$, it would mean that there are two different unit circles with $3$ different intersection points, which is impossible.
\end{pf}

\begin{rem}[description={Globus, Parshall, 2019}]
{\rm{\cite{gp}}} A complete classification of the minimal forbidden graphs (graphs which are not unit distance graphs, but all their real subgraphs are) on at most $9$ vertices was done.
\end{rem}

\begin{thm}[description={Schaefer, 2013}]
{\rm{\cite{schae}}} To decide, whether a graph is a UDG or not is NP-hard, more specifically ETR-complete.
\end{thm}

\begin{defn}
Let us call the graph on Figure \ref{moser} a \nt{Moser spindle}.
\end{defn}

\begin{prop}
The Moser spindle is a faithful unit distance graph.
\end{prop}

\begin{pf}

\begin{center}
	\input{Figures/Moserspindle}
	\captionof{figure}{}
	\label{moser}
\end{center}

Let us construct the point set in Figure \ref{moser}, and we will show that it is a FUDR of a Moser spindle. Let us take an equilateral triangle $ABC$ and another one, $CBD$. If we rotate $B$, $C$ and $D$ around $A$ with an angle $2\cdot\arcsin\left({1\over2\sqrt{3}}\right)$, we get $E$, $F$ and $G$, respectively and the point set $\lbrace A,E,F,G\rbrace$ will be isometric to the point set $\lbrace A,B,C,D\rbrace$, but also the isosceles triangle $ADG$ has a base length $2\cdot{\sin{\arcsin\left({1\over2\sqrt{3}}\right)}}\cdot\sqrt{3}=2\cdot{1\over2\sqrt{3}}\cdot\sqrt{3}=1$. So the connected vertices indeed have distance $1$. We now only have to prove that the remaining pairs have some distance different from $0$ and $1$. Sets $\lbrace A\rbrace$, $\lbrace B,E,C,F\rbrace$ and $\lbrace D,G\rbrace$ are pairwise disjoint as $A$ has distance $0$ from itself, and all of the points in $\lbrace B,E,C,F\rbrace$ have distance $1$ from $A$, while those in $\lbrace D,G\rbrace$ have distance $\sqrt{3}$ from $A$. And there are no other $1$ distances between these sets than what we already know since the only possible pairs would be $BG$, $ED$, $CG$ and $FD$, but they would make the induced subgraphs of $\lbrace A,G,B,E,F\rbrace$, $\lbrace A,D,B,E,C\rbrace$, $\lbrace A,G,B,C,F\rbrace$ and $\lbrace A,D,B,C,F\rbrace$ $K_{2,3}$'s respectively, which is impossible. So we now only have to examine the possible pairs inside sets $\lbrace A\rbrace$, $\lbrace B,E,C,F\rbrace$ and $\lbrace D,G\rbrace$, but among these only $\lbrace B,E,C,F\rbrace$ contains non-adjacent pairs: $BE$, $CF$ and $BF$ and since they are all chords on the unit circle around $A$ and their central angles are $\angle{BAE}=\angle{CAF}=60^\circ-2\arcsin\left({1\over2\sqrt{3}}\right)$ and $\angle{BAF}=60^\circ+2\arcsin\left({1\over2\sqrt{3}}\right)$, so their length cannot be $0$, neither $1$. So the described (and drawn) point set indeed is a FUDR of the Moser spindle.
\end{pf}

\begin{defn}
We call a unit distance graph \nt{weakly rigid}, if it has only one unit distance representation up to isometry.
\end{defn}

\begin{defn}
We call a unit distance graph \nt{strongly rigid}, if it has only one unit distance representation up to isometry, even if we label its points.
\end{defn}

\begin{rem}
These definitions do not fully correspond with those in rigidity theory.
\end{rem}

\begin{prop}[label=merev3szogracs]
Suppose that a graph $G$ is the union of triangles so that one can walk from any triangle to any other by taking steps only between triangles with a common side. If such a graph is a unit distance graph, it is strongly rigid.
\end{prop}

\begin{pf}
It can be seen by induction:

(If the graph is the union of $0$ triangles, it does not have any vertices or edges.)

For $1$ triangle, it is trivial: the graph can only be a $K_3$ and all equilateral triangles are isometric with each other.

If we have seen it for $n-1$ triangles for some $n\ge 2$, it is also true for a graph $G$, which is the union of $n$ path-connected triangles: we can find a spanning tree in the graph of triangles, where the edges are the adjacencies defined above. And if we take a leaf $L$ in this graph, then it is possible that all of its vertices are among the vertices of the other $n-1$ triangles, in which case, the graph is the union of $n-1$ triangles, which are also connected by paths because deleting a leaf from the spanning tree still leaves it connected. The other possibility is that $L$ has exactly two common vertices with the other $n-1$ triangles: the endpoints of the common side with its neighbouring triangle. But in this case, one equilateral triangle having these two points as its side is already in the point set, so we can place the remaining vertex of $L$ only at one place. So the induction step is true in this case too.
\end{pf}

\begin{prop}
The Moser spindle is weakly rigid, but it is not strongly rigid.
\end{prop}

\begin{pf}
Graphs $ABCD$ and $AEFG$ on Figure \ref{moser} are strongly rigid because of Proposition \thmref{merev3szogracs}, and the angle in which they are rotated to each other is also determined, but the orientation of the rotation is not. Thus, the point set is determined, but Figure \ref{moser} and Figure \ref{moser2} are examples of two different labellings.
\end{pf}

\begin{center}
	\input{Figures/Moserspindle2}
	\captionof{figure}{}
	\label{moser2}
\end{center}

\begin{prop}[label=maxelszam5ig]
The maximal edge numbers of a unit distance graph with $1$, $2$, $3$, $4$ and $5$ vertices are $0$, $1$, $3$, $5$, $7$.
\end{prop}

\begin{pf}
For $1$, $2$ and $3$ vertices these are the maximal edge counts and it can be trivially obtained by taking a point, a segment of length $1$ and an equilateral triangle of side-length $1$, respectively.

For $4$ vertices, two equilateral triangles with side-length $1$ and a common side give an example for $5$ edges (Figure \ref{4csucs5el}). But there cannot be $6$ edges, because of Proposition \thmref{k4}

For $5$ vertices and $7$ edges, the point set on Figure \ref{5csucs7el} is a good example. But if a unit distance graph with $5$ vertices and at least $8$ edges existed, then the sum of the degrees of the vertices would be at least $16$, which would mean that there would be at least one vertex, which has degree $4$ (otherwise the sum of the degrees would be at most $5\cdot 3=15$), let it and its image be called $A$. Then all images of all the other vertices would be on a circle around $A$, and also the sum of their degrees would be at least $16-4=12$, so either all four of them have degree $3$ or at least one of them (say, $B$) also has degree $4$. But the latter would mean that the graph contains a $K_{2,3}$ with $A$ and $B$ on one side and $C$, $D$ and $E$ on the other, which is impossible because of Proposition \thmref{k23}. So the former case is true, but then all of $B$, $C$, $D$ and $E$ have their rotated image with both $60^\circ$ and $-60^\circ$ around $A$ among the points. But this would mean that starting from $B$, we would draw Figure \ref{fudgreszenemfudg}, which has $7$ vertices, so it is a contradiction.
\end{pf}

\begin{center}
	\input{Figures/maximalis5csucsuUDG}
	\captionof{figure}{}
	\label{5csucs7el}
\end{center}

\begin{defn}
Let us call a planar point set a \nt{planar unit-hypercube of order $n$} for some $n$ non-negative integer if it has $2^n$ distinct points and taking one of the points as the origin, the points can be written in a form $\sum_{i=1}^n{\lambda_i\cdot v_i}$, where $v_i\in\mathbb R^2, \abs{v_i}=1$ $(i=1,...,n)$ and $\lambda_i\in\lbrace 0,1\rbrace$.
\end{defn}

\begin{thm}
Planar unit-hypercubes exist for any order $n$.
\end{thm}

\begin{pf}
If we have created a planar unit-hypercube for $n-1$ (it is trivial for $0$), we can take the same vectors for $v_1, v_2, v_3, ..., v_{n-1}$ and find a $v_n$ such that none of the two sub-hypercubes generated by $v_n$ (the set of points for which $\lambda_n=0$ and for which $\lambda_n=1$) have a common point. We can do it, because every pair of points from the two different sub-hypercubes can coincide in at most one case out of the continuously many. So by induction, we can find a planar unit-hypercube of order $n$ for all $n$'s.
\end{pf}

\begin{thm}
The edge-graph of an $n$-dimensional hypercube is a unit distance graph for all $n\in\mathbb N$ and the planar unit-hypercubes of order $n$ are its only unit distance representations.
\end{thm}

\begin{pf}
A planar unit-hypercube $P$ of order $n$ is a UDR of an $n$-dimensional hypercube, since if we have a representation as described above, then we can make a bijection with the vertices of the unit cube $\left[0,1\right]^n$ such that for every point in $S$, we take the $\lambda_i$'s as $i$th coordinates for all $1\le i\le n$. This makes a bijection as every combination of the possible $\lambda$'s is ocurring among the points of $S$ (otherwise it could not have $2^n$ points) and the vertices of the unit hypercube are also the vertices with these coordinates. And those points of $S$, whose representation only differs in $\lambda_j$ for one $j$ have a distance of $1$, since that is the length of their difference vector $v_j$. And this is exactly the case, when their counterparts in the hypercube have an edge between them, so it is true that for every edge, the corresponding points in $P$ have distance $1$. And since planar unit-hypercubes exist for any order $n$, it means that the $n$-dimensional hypercubes are UDRs.

Also, if a hypercube has a UDR, let us choose an arbitrary vertex $A$ of it and let its image be the origin and let the $v_i$ vectors be the vectors from the origin to its neighbours.  Now we only have to prove that the only way to draw an UDR for a hypercube is to draw every vertex at the sum of a subset of the $v_i$ vectors, which is true by induction: If for all vertices with a graph distance at most $n-1$ from the origin, we have found a sum of at most $n-1$ of the $v_i$'s, so that this makes a UDG, we have to prove that it also works for $n$ (the starting step is trivial: for $n=1$ it is true). And the induction step is true: for two vertices $P$ and $Q$, which have graph distance $n-1$ from $A$ and graph distance $2$ from each other, we have already drawn their neighbour $R$ closer to $A$, so the farther one, $S$ only can be at the remaining vertex of the rhombus generated by $P$, $R$ and $Q$ and since $\vec{RP}$ and $\vec{RQ}$ are from the $v_i$'s and all of them have a $\lambda$-representation, $S$ also has one.
\end{pf}

\begin{thm}[description={Schade, 1993}]
{\rm{\cite{sch}}} The maximum edge numbers of a unit distance graph with $n$ vertices for $n\le14$ are the following:
\end{thm}

\begin{center}
\begin{tabular}{ | c | c | c | }
\hline
number of points&maximal edge count&example\\
\hline
1&0&\usetikzlibrary{arrows}
\begin{tikzpicture}[line cap=round,line join=round,scale=0.7]
\clip(-1.5,-0.5) rectangle (-0.5,0.5);
\begin{scriptsize}
\fill [color=black] (-1,0) circle (1.5pt);
\end{scriptsize}
\end{tikzpicture}
\\
\hline
2&1&\usetikzlibrary{arrows}
\begin{tikzpicture}[line cap=round,line join=round,scale=0.7]
\clip(-1.25,-0.25) rectangle (0.25,0.25);
\draw (-1,0)-- (0,0);
\begin{scriptsize}
\fill [color=black] (0,0) circle (1.5pt);
\fill [color=black] (-1,0) circle (1.5pt);
\end{scriptsize}
\end{tikzpicture}
\\
\hline
3&3&\usetikzlibrary{arrows}
\begin{tikzpicture}[line cap=round,line join=round,scale=0.7]
\clip(-1.25,-0.25) rectangle (0.25,1.1);
\draw (-0.5,0.87)-- (-1,0);
\draw (-0.5,0.87)-- (0,0);
\draw (-1,0)-- (0,0);
\begin{scriptsize}
\fill [color=black] (0,0) circle (1.5pt);
\fill [color=black] (-1,0) circle (1.5pt);
\fill [color=black] (-0.5,0.87) circle (1.5pt);
\end{scriptsize}
\end{tikzpicture}
\\
\hline
4&5&\input{Figures/Schade4}\\
\hline
5&7&\input{Figures/Schade5}\\
\hline
6&9&\input{Figures/Schade6-1}\input{Figures/Schade6-2}\input{Figures/Schade6-3}\input{Figures/Schade6-4}\\
\hline
7&12&\input{Figures/Schade7}\\
\hline
8&14&\input{Figures/Schade8-1}\input{Figures/Schade8-2}\input{Figures/Schade8-3}\\
\hline
9&18&\input{Figures/Schade9}\\
\hline
10&20&\input{Figures/Schade10}\\
\hline
11&23&\input{Figures/Schade11-1}\input{Figures/Schade11-2}\\
\hline
12&27&\input{Figures/Schade12}\\
\hline
13&30&\input{Figures/Schade13}\\
\hline
14&33&\input{Figures/Schade14-1}\input{Figures/Schade14-2}\\
\hline
\end{tabular}
\end{center}

The graphs above are the only optimal examples for $n\le 13$ and two optimal examples for $n=14$.

The rest of the proof goes through all the possible cases with the help of Proposition \thmref{udgedgedensity}, Proposition \thmref{minkowski} and an estimate using the minimal degree vertex of the maximal UDGs, among others.

\begin{center}
\begin{minipage}{.49\textwidth}
	\centering
	\input{Figures/Schade15}
	\captionof{figure}{}
	\label{schade15}
\end{minipage}
\begin{minipage}{.49\textwidth}
	\centering
	\input{Figures/Schade16}
	\captionof{figure}{}
	\label{schade16}
\end{minipage}%
\end{center}

The above are the maximal known UDGs for vertex number $15$ and $16$. The latter is a planar unit-hypercube of order $4$.

Let us call the graphs in the table Graph $n$.$k$, where $n$ denotes the number of vertices and $k$ denotes their place in the cell in the table they are in from left to right in case there are more than $1$ graphs for the given $n$. In case there is only one graph listed for some vertex number $n$, let us call it Graph $n$. Let us also call the above two graphs Graph 15 and Graph 16, respectively.

\begin{thm}
The above graphs are unit distance graphs.
\end{thm}

\begin{pf}
Here are the coordinates of Graph 16:

\begin{center}
	\input{Figures/16csucsuegyseggraf}
	\captionof{figure}{}
	\label{16csucs41el}
\end{center}

$A=\left(-{1\over2},0\right)$

$B=\left({1\over2};0\right)$

$C=\left(0,\sqrt{11\over4}\right)$

$D=\left(0,-\sqrt{11\over4}\right)$

$E=\left(-{1\over4}-\sqrt{11\over48},-\sqrt{11\over16}-\sqrt{1\over48}\right)$

$F=\left({1\over4}-\sqrt{11\over48},-\sqrt{11\over16}+\sqrt{1\over48}\right)$

$G=\left(-{1\over4}+\sqrt{11\over48},-\sqrt{11\over16}+\sqrt{1\over48}\right)$

$H=\left({1\over4}+\sqrt{11\over48},-\sqrt{11\over16}-\sqrt{1\over48}\right)$

$I=\left({1\over4}+\sqrt{11\over48},\sqrt{11\over16}+\sqrt{1\over48}\right)$

$J=\left(-{1\over4}+\sqrt{11\over48},\sqrt{11\over16}-\sqrt{1\over48}\right)$

$K=\left({1\over4}-\sqrt{11\over48},\sqrt{11\over16}-\sqrt{1\over48}\right)$

$L=\left(-{1\over4}-\sqrt{11\over48},\sqrt{11\over16}+\sqrt{1\over48}\right)$

$M=\left(-\sqrt{11\over12},0\right)$

$N=\left(\sqrt{11\over12},0\right)$

$O=\left(0,-\sqrt{1\over12}\right)$

$P=\left(0,\sqrt{1\over12}\right)$

Segment $AB$ trivially has length $1$.

Segments $DE$, $FM$, $GA$, $HO$, $PL$, $BK$, $NJ$ and $IC$ have length 
$$\sqrt{\left(-{1\over4}-\sqrt{11\over48}\right)^2+\left(\sqrt{11\over16}-\sqrt{1\over48}\right)^2}=\sqrt{{1\over16}+\sqrt{11\over192}+{11\over48}+{11\over16}-\sqrt{11\over192}+{1\over48}}=\sqrt{48\over48}=1.$$

Segments $DH$, $GN$, $FB$, $EO$, $PI$, $AJ$, $MK$ and $LC$ have length
$$\sqrt{\left({1\over4}+\sqrt{11\over48}\right)^2+\left(\sqrt{11\over16}-\sqrt{1\over48}\right)^2},$$which is similarly $1$.

Segments $DF$, $HB$, $AL$, $OK$, $EM$, $GP$, $JC$ and $NI$ have length
$$\sqrt{\left({1\over4}-\sqrt{11\over48}\right)^2+\left(\sqrt{11\over16}+\sqrt{1\over48}\right)^2}=\sqrt{{1\over16}-\sqrt{11\over192}+{11\over48}+{11\over16}+\sqrt{11\over192}+{1\over48}}=\sqrt{48\over48}=1.$$

Segments $DG$, $HN$, $EA$, $ML$, $FP$, $OJ$, $BI$, $KC$ and $ML$ have length
$$\sqrt{\left(-{1\over4}+\sqrt{11\over48}\right)^2+\left(\sqrt{11\over16}+\sqrt{1\over48}\right)^2},$$
which is similarly $1$.

Segments $EG$, $MP$, $ON$ and $KI$ have length $\sqrt{\left(\sqrt{11\over12}\right)^2+\left(\sqrt{1\over12}\right)^2}=\sqrt{{11\over12}+{1\over12}}=\sqrt{12\over12}=1$.

Segments $FH$, $MO$, $PN$ and $LJ$ have length $\sqrt{\left(\sqrt{11\over12}\right)^2+\left(-\sqrt{1\over12}\right)^2}$, which is similarly $1$.

So all the segment lengths are $1$.

From the above, it also can be seen that Graph 16 is a planar unit-hypercube, which has 9 non-hypercube edges: its edge vectors are $\vec{DE}$, $\vec{DF}$, $\vec{DG}$ and $\vec{DH}$.

And all of the graphs listed above for $9\le n\le 16$, Graph 6.4, Graph 8.2 and Graph 8.3 are subgraphs of this graph:

Graph 6.4 is the subgraph induced by $\lbrace C,L,J,K,M,D\rbrace$.

Graph 8.2 is induced by $\lbrace C,L,J,K,M,	A,O,E\rbrace$.

Graph 8.3 is the subgraph induced by $\lbrace M,O,N,F,G,E,H,D \rbrace$.

Graph 9 is the subgraph induced by $\lbrace C,L,I,K,J,P,M,N,O\rbrace$.

Graph 10 is the subgraph induced by $\lbrace C,L,I,K,J,P,M,N,O,B\rbrace$.

Graph 11.1 is the subgraph induced by $\lbrace C,L,I,K,J,P,M,N,O,B,\allowbreak F\rbrace$.

Graph 11.2 is the subgraph induced by $\lbrace C,L,I,K,J,P,M,N,O,B,A\rbrace$.

Graph 12 is the subgraph induced by $\lbrace C,L,I,K,J,P,M,N,O,B,\allowbreak F,H\rbrace$.

Graph 13 is the subgraph induced by $\lbrace C,L,I,K,J,P,M,N,O,B,\allowbreak F,H,A\rbrace$.

Graph 14.1 is the subgraph induced by $\lbrace L,I,K,J,P,M,N,O,B,\allowbreak F,H,A,G,E\rbrace$.

Graph 14.2 is the subgraph induced by $\lbrace C,L,I,K,J,P,M,N,O,B,F,H,A,G\rbrace$.

So these are also unit distance graphs.

All the other graphs in the table above are triangular grid graphs, so they also are unit distance graphs.
\end{pf}

\begin{thm}[label=schaderigid]
Graphs 1, 2, 3, 4, 5, 6.1, 6.2, 6.3, 7 and 8.1 are strongly rigid, while Graphs 6.1, 8.2, 8.3, 9 and 10 are not even weakly rigid.
\end{thm}

\begin{pf}
The triangular grids are strongly rigid because of Proposition \thmref{merev3szogracs}

Graph 6.4 is not rigid as triangles $CLJ$ and $KMD$ are rigid, but their difference vector can be any vector of length $1$, with the exception of the side vectors of the triangles themselves.

Graph 8.2 is not rigid either as subgraphs $CLJA$ and $KMOE$ are rigid themselves but can move independently from each other so that they are each other's translated images for some vector of length $1$ again not one of the vectors among their own vertices.

Graph 8.3, Graph 9 and Graph 10 are not rigid either:

If we choose an arbitrary rhombus as $MOJL$, then by putting equilateral triangles on its sides as it is seen in the drawing in the table, we create $N$, $P$, $K$ and $C$. Since the first two vertices of triangle $LMP$ and $JON$ have the same vector ($\vec{MO}=\vec{LJ}$) as their difference (because of $MOJL$ being a rhombus), it is true for the third vertices $P$ and $N$ too (so $\vec{PN}=\vec{MO}=\vec{LJ}$). Similarly $\vec{KC}=\vec{ML}=\vec{OJ}$. So $PN$ and $KC$ have distance one, so all edges have length $1$ and there are infinitely many ways to choose $MOJL$ so that no points coincide, as a sufficiently small change to the drawn rhombus does not make any two points coincide. So Graph 8.3 is not rigid.

Let us now draw two equilateral triangles to $PN$ and $CK$ as in the drawing and let us name their third vertices $I_1$ and $I_2$, respectively. Since triangle $PNI_1$ is a translated image of $MOK$ with $\vec{MP}$, while triangle $CKI_2$ is also a translated image of $LMP$ with $\vec{MK}$, so $\vec{I_1I_2}=\vec{IP}+\vec{PM}+\vec{MK}+\vec{KI}=\vec{KM}+\vec{PM}+\vec{MK}+\vec{MP}=0$, so $I_1=I_2$ regardless of how we chose $MOJL$. And again, if we change $\angle{OML}$ with a sufficiently small angle, none of the points will coincide. This makes Graph 9 not rigid too.

And if we draw triangle $IKB$ to side $IK$, we get $B$, which also does not coincide with any other points if we change $\angle{OML}$ with a sufficiently small angle, so Graph 10 is not rigid too.

It is also worth mentioning, why Graph 9 has $4$ more edges than Graph 8, while Graph 10 has only $2$ more edges than Graph 9, and also why the 10th vertex in Graph 10 can seemingly put in only one place. Graph 9 is strongly regular with parameters $(9,4,1,2)$, and also it is vertex-transitive and edge-transitive and it is a generalized quadrangle known as $GQ(2,1)$. So first, this means that all the subgraphs of Graph 9 with $8$ vertices are isometric to Graph 8 (as we know that it is a subgraph), so the inequality from Proposition \thmref{udgedgedensity} holds with equality, which not only explains the fact that it has $4$ more edges than Graph 8, but also makes it understandable why we cannot add more edges: Graph 9 is a very "perfectly constructed" graph, so it is more profitable to add a vertex to it even if the maximal degree we can make it to have is $2$, than constructing another graph that does not have Graph 9 as a subgraph. But it still would be possible that we can add a vertex with degree more than $2$ to Graph 9 as it is not rigid, moreover it has continuously many UDRs. But this is not possible: if the new vertex would have two neighbours, which are not adjacent, we would get $K_{2,3}$ because all of the non-adjacent vertices of Graph 9 have two common neighbours because of the fourth parameter of the strong regularity. So all of the neighbours of the new vertex are adjacent, which means that if it would have degree more than $2$, we would get a $K_4$ in the graph, which is impossible. And this also explains why is Graph 10 the only maximal UDG for $10$ vertices: we only can add a vertex with degree $2$ to the endpoints of an edge and Graph 9 is edge-transitive so all the graphs we get this way are isomorphic to each other, although their UDR is not fixed even if we fix the UDR of Graph 9.
\end{pf}

\begin{thm}[label=legfeljebb5csucsugrafok]
The non-empty unit distance graphs with at least $1$, but at most $5$ vertices are the following (they are also faithful unit distance graphs):

\begin{center}
\begin{minipage}{\textwidth}
	\input{Figures/Egyseggrafok}
	\captionof{figure}{}
	\label{egyseggrafok}
\end{minipage}
\end{center}

(The above drawings are FUDRs of the graphs, so we do not need any further proof that they are unit distance graphs.)
\end{thm}

\begin{pf}
The cases for $1$, $2$ and $3$ vertices (Graphs 1-7) are true because all such graphs are unit distance graphs.

In case of $4$ vertices: 

For $0$ edges, Graph $8$ is the only possibility and for $1$ edge, Graph $9$ is the only one. For $2$ edges, either there exists a vertex with degree $2$, in which case we get Graph $10$ or there does not exist such a vertex, in which case we get Graph $11$. For $3$ edges, either there exists a vertex with degree $3$, in which case we get Graph $12$, or the maximum degree is $2$ and there are $3$ such vertices, in which case we get Graph $13$, or the maximum degree is $2$ and there are $2$ such vertices, in which case we get Graph $14$. If there would be at most $1$ vertex with degree $2$, then the sum of the degrees would be at most $2+3\cdot1=5$, which is a contradiction. For $4$ edges, we get the complement of a graph with $2$ edges, so either we get Graph $15$ (the complement of Graph $10$) or Graph $16$ (the complement of Graph $11$). For $5$ edges, we only can get Graph $17$, as it is the complement of Graph $9$. For $6$ edges, there are no unit distance graphs, since the only possibility would be $K_4$, which was already shown not to be a unit distance graph in Proposition \thmref{k4}.

In case of $5$ vertices:

For $0$ edges, the only possible graph is Graph $18$, and for $1$ edges the only possible graph is Graph $19$. For $2$ edges, either there is a vertex with degree $2$, in which case we get Graph $20$, or there is not, in which case we get Graph $21$.

For $3$ edges, if the maximal degree is $3$, we get Graph $22$. If it is $2$ and there are $3$ vertices with degree $2$, then we get Graph $23$, if there are $2$ such vertices, we get Graph $24$ and if there is only $1$, then we get Graph $25$. There cannot be less than $1$ or more than $3$ such vertices, as the sum of the degrees is $6$.

For $4$ edges, if the maximal degree is $4$, we get Graph $26$, if it is $3$, then we either get Graph $27$ or Graph $28$ depending on the position of the remaining edge, if it is $2$, then there are at most $4$ vertices with degree $2$. If there are exactly $4$, the only possibility is Graph $29$, and if there are $3$, we either get Graph $30$ or Graph $31$. And since having at most $2$ vertices with degree $2$ and none with more would make the sum of the degrees at most $2\cdot2+3\cdot1=7$, this would be a contradiction, so we do not have any more cases.

For $5$ edges, either there is a vertex with degree $4$ and we get Graph $32$, or the maximal degree is $3$ in which case we get Graph $33$, Graph $34$, Graph $35$ or Graph $36$ depending on the position of the other two edges, or all of the vertices have degree at most $2$. But since the sum of the degrees is $10$, the only such possibility is that all $5$ vertices have degree $2$, which implies that the graph is $C_5$, that is, Graph $37$.

For $6$ edges, we have to go through the complements of the graphs with $4$ edges (for $4$ edges, all of the graphs were unit distance graphs, so we will not miss anything). The complement of Graph $26$ is $K_4$, which is not a unit distance graph as it was proved in Proposition \thmref{k4}. The complement of Graph $27$ is Graph $38$. The complement of Graph $28$ is Graph $39$. The complement of Graph $29$ is Graph $40$. The complement of Graph $30$ is Graph $41$. The complement of Graph $31$ is $K_{2,3}$, which is not a unit distance graph by Proposition \thmref{k23}.

For $7$ edges, we have to go through the complements of the graphs with $3$ edges. The complement of Graph $22$ contains a $K_4$, so it is not a unit distance graph. The complement of Graph $23$ contains a $K_{2,3}$, so it is not a unit distance graph. The complement of Graph $24$ is Graph $42$. The complement of Graph $25$ is the graph contains a $K_{2,3}$, so it is not a unit distance graph.
\end{pf}

\section{Specific kinds of unit distance graphs}

\begin{defn}
Let us take a unit distance graph, for which there exists a unit distance representation in which the edges do not cross each other. Let us call such a graph a \nt{matchstick graph.}
\end{defn}

\begin{defn}
Let us take a unit distance graph, which has a unit distance representation with its vertices being the vertices of a convex polygon. Let us call such a graph a \nt{convex unit distance graph} (CUDG) and any such representation a convex unit distance representation (CUDR).
\end{defn}

\begin{defn}
Let us call the maximal possible edge count of a convex unit distance graph called $u_c(n)$.
\end{defn}

\begin{thm}[description={Füredi, 1990}] {\rm{\cite{f}}} $u_c(n)=O(n\log{n})$.
\end{thm}

\begin{pf}
The following proof comes from Peter Braß and János Pach{\rm{\cite{bp}}}:

Let $G$ be a CUDG with $n$ vertices and $u_c(n)$ edges and let us take a $P$ point set, which is a CUDR for $G$ and whose vertices are $p_1,p_2,p_3,...,p_n$, where the ordering is defined by their location on their convex hull. Also, let us choose $p_i$ as an antipodal point of $p_1$. Now we can suppose that $i>{n\over2}$, otherwise we can reverse the ordering of the points. Let us also take parallel lines $l$ and $l"$, which go through $p_1$ and $p_i$ and all the other points are inside the open strip between them. We can take them so that no segment between the points of $P$ are parallel with $l$ and $l'$, and no segment between the points of $P$ is perpendicular to them. Let us now define a coordinate system so that $l$ and $l'$ are parallel with the $x$ axis so that $p_1$ is the lowest point of $P$. We now can divide the edges of $G$ into two parts: let those be red edges, whose image is a segment with a positive slope and let those be blue edges, whose image is a segment with a negative slope. Let us now take those edges, whose image does not cross $p_1 p_i$ and are red. If we assign all such edges on the right of $p_1 p_i$ to their right endpoint, while those on the left of $p_1p_i$ to their left endpoint, then all points have at most one edge assigned to them: we can suppose WLOG that the points were numbered counterclockwise, which means that for some $1\le j\le i$, all points $p_k$ from $P$ that have a smaller $x$ and a smaller $y$ coordinate (that is what $p_j$ being the right endpoint of a red edge connecting it with that certain point means) have a smaller index than $j$ and all the possible $p_k$'s are on the part of the broken line $p_1p_2...p_i$, where both coordinates are increasing, so their distance from $p_j$ is strictly decreasing meaning that at most one of them having a distance of $1$ from $p_j$. The same applies on the other side of the line $p_1p_i$, so indeed there are at most $n$ red edges not crossing $p_1p_i$, and the same also applies for blue edges, so in total at most $2n$ edges do not cross line $p_1p_i$.

And similarly, this is also true for line $p_{\lceil{i\over2}\rceil} p_l$, where $p_l$ is an antipodal point of $p_{\lceil{i\over2}\rceil}$. And the two lines divide $P$ into four (not necessarily non-empty) parts, so apart from at most $4n$ edges, all edges cross both meaning that defining $P_1$ as $\lbrace p_2,p_3,...,p_{\lceil{i\over2}\rceil-1}\rbrace\cup\lbrace p_{i+1},p_{i+2},...,p_{j-1}\rbrace$ and $P_2$ as $\lbrace p_{\lceil{i\over2}\rceil+1},p_{\lceil{i\over2}\rceil+2},...,p_{i-1}\rbrace\cup\lbrace p_{j+1},p_{j+2},...,p_n\rbrace$, all edges which cross both lines have their endpoints either both from $P_1$ or both from $P_2$. And $\abs{P_1}+\abs{P_2}$ ($P_1$ and $P_2$ contain all points of $P$ exactly once except for $p_1$, $p_{\lceil{i\over2}\rceil}$, $p_i$ and $p_l$, which are not contained in any of them) and $\min{\left(\abs{P_1},\abs{P_2}\right)}\ge{n-7\over4}$. From this, we can conclude the statement of the theorem by induction.
\end{pf}

\chapter{$(1,d)$-graphs}

In this chapter, we will discuss several statements about $(1,d)$-graphs, most of them are fuelled by the probabilistic approach of the Hadwiger--Nelson problem to be discussed in Chapter~4.

\begin{defn}
A coloured graph with $2$ colours is called a \nt{$(1,d)$-graph} for some $d$, if it can be represented by a point set in the plane for which those points have distance $1$, which have an edge coloured with colour $1$ between them, and those have distance $d$, which have an edge of colour $2$ between them. From now on let us call these edges $1$-edges and $d$-edges.
\end{defn}

\begin{defn}
An uncoloured graph is called a \nt{$(1,d)$-graph} for some $d$ if it has a colouring that is a $(1,d)$-graph for that certain $d$.
\end{defn}

\begin{defn}
A coloured graph with $2$ colours is called a \nt{faithful $(1,d)$-graph} for some $d$, if it can be represented by a point set in the plane for which exactly those points have distance $1$, which have an edge coloured with colour $1$ between them, and exactly those are of distance $d$, which have an edge of colour $2$ between them. We also call these edges $1$-edges and $d$-edges.
\end{defn}

\begin{defn}
An uncoloured graph is called a \nt{faithful $(1,d)$-graph} for some $d$ if it has a colouring that is a faithful $(1,d)$-graph for some $d$.
\end{defn}

\begin{defn}
Let us take a bicoloured graph $G$ and let us call the set of the numbers $d$ for which $G$ is a $(1,d)$-graph the \nt{range} of the graph.
\end{defn}

\begin{defn}
The \nt{inverse} of a $(1,d)$-graph is a $(1,d')$-graph for some $d'$ which is isomorphic to the original one, but has its $1$-edges and $d$-edges switched.
\end{defn}

\begin{rem}
The inverse of the inverse of a graph $G$ is $G$ itself, because it will have the same edges, and instead of the $1$-edges, there will be $1$-edges and instead of the $d$-edges, there will be $d$-edges.
\end{rem}

Note that if $G$ is a $(1,d)$-graph, then its inverse is a $(1,d')$-graph for $d'={1\over d}$.

\begin{thm}
The numbers $d$ for which a two-coloured graph is $(1,d)$-graph form a finite number of closed intervals.
\end{thm}

\begin{pf}
The good $d$'s can be written as the solutions of a polynomial equation and thus, by the Tarski--Seidenberg theorem {\rm{\cite{ts}}} the statement concludes.
\end{pf}

\begin{defn}
Let us call a bicoloured graph with a range that includes intervals of positive length a \nt{$(1,d)$-graph for a general $d$}.
\end{defn}

\begin{ex}[label=nehany1dgraf]
Here are a few examples for $(1,d)$-graphs with $5$ vertices and $8$ edges:
\end{ex}

\begin{center}
\begin{minipage}{.24\textwidth}
	\centering
	\input{Figures/1dgraph1}
	\captionof{figure}{$d\in\left(0,2\right]\setminus\lbrace1\rbrace$}
	\label{1dgraph1}
\end{minipage}%
\begin{minipage}{.24\textwidth}
	\centering
	\input{Figures/1dgraph2}
	\captionof{figure}{$d\in\left[{1\over2},\infty\right)\setminus\lbrace1\rbrace$}
	\label{1dgraph2}
\end{minipage}
\begin{minipage}{.24\textwidth}
	\centering
	\input{Figures/1dgraph3}
	\captionof{figure}{$d\in\left[{1\over2},\infty\right)\setminus\lbrace1\rbrace$}
	\label{1dgraph3}
\end{minipage}
\begin{minipage}{.24\textwidth}
	\centering
	\input{Figures/1dgraph4}
	\captionof{figure}{$d\in\left(0,2\right]\setminus\lbrace1\rbrace$}
	\label{1dgraph4}
\end{minipage}
\end{center}
\begin{center}
\begin{minipage}{.24\textwidth}
	\centering
	\input{Figures/1dgraph5}
	\captionof{figure}{$d\in\left[{1\over2},2\right]\setminus\lbrace1\rbrace$}
	\label{1dgraph5}
\end{minipage}%
\begin{minipage}{.24\textwidth}
	\centering
	\input{Figures/1dgraph6}
	\captionof{figure}{$d={\sqrt{6}+\sqrt{2}\over2}$}
	\label{1dgraph6}
\end{minipage}
\begin{minipage}{.24\textwidth}
	\centering
	\input{Figures/1dgraph7}
	\captionof{figure}{$d={1+\sqrt{5}\over2}$}
	\label{1dgraph7}
\end{minipage}
\begin{minipage}{.24\textwidth}
	\centering
	\input{Figures/1dgraph8}
	\captionof{figure}{$d={1-\sqrt{5}\over2}$}
	\label{1dgraph8}
\end{minipage}
\end{center}

Red edges mean $1$-edges and blue edges mean $d$-edges.

\begin{prop}[label=1dsokszog]
If a $(1,d)$-graph contains an $n$-cycle for some $n$ with only one $1$-edge, then its range does not contain any $d<{1\over n-1}$, and if it contains an $n$-cycle with only one $d$-edge, then its range does not contain any $d>n-1$.
\end{prop}

\begin{pf}
For the first part: If there would be a $(1,d)$-representation of such a $(1,d)$-graph for some $d<{1\over n-1}$, that would mean that there is a broken line formed by $n-1$ segments of length smaller than ${1\over n-1}$ and whose starting point and endpoint have distance $1$, which is impossible because of the triangle inequality. For the second part, if such a graph would be a $(1,d)$-graph for some $d>n$, then its inverse graph would be a $(1,d')$-graph for $d'={1\over d}<n-1$.
\end{pf}

\begin{prop}[label=1dgraf58]
A bicoloured graph with $5$ vertices and at least $8$ edges cannot be a $(1,d_1)$-graph for some $d_1<{1\over2}$ and a $(1,d_2)$-graph for some $d_2>2$ at the same time.
\end{prop}

\begin{pf}
Let us make pairs so that the first element of the pair is a $3$-element set of vertices of the graph, while the other is an edge running between these three vertices. There are at least $3\cdot8=24$ such pairs, because for all edges, we can find a third vertex in $3$ possible ways. So there are at least $4$ triangles in the graph, since for $3$ element vertex sets, we have got $3$ such pairs if they form a $K_3$ in the graph and we have got $2$ such pairs if they do not, so in total, we would have at most $3\cdot3+2\cdot7=23$ such pairs without $4$ triangles.

\begin{lem}[label=szabalyos3szogek]
There cannot be $4$ monochromatic triangles among $5$ points in the representation of a $(1,d)$-graph (or equivalently, there cannot be $4$ equilateral triangles with side length $1$ or $d$ in such a graph, as we can assume the representation is faithful).
\end{lem}

\begin{pf}
Let us suppose that there are at least $3$ 

of them, which have side length $1$.

In this case, these $3$ triangles have $9$ vertices in total (counted with multiplicity), so if for all three pairs of them, there would only be at most one common vertex, then there would be at least $6$ points in total. So there are two triangles out of these three, which have exactly two common vertices, so they have one common side (there cannot be three common vertices, because, then these two would be the same triangle). Let these be $ABC$ and $BCD$. In this case, the third triangle has at least $2$ vertices from $\lbrace A,B,C,D\rbrace$, since there is only one point left. Out of $A$, $B$, $C$ and $D$, only $A$ and $D$ do not have a distance that is not $1$, but $BC$ does not have any more equilateral triangles having it as a side other than $ABC$ and $BCD$, so it is not also possible that these are the two points. So we have got four possibilities for the two common points of the third triangle with $\lbrace A,B,C,D\rbrace$: $(A,B)$, $(B,D)$, $(D,C)$ and $(C,B)$. But because of the symmetry of the rhombus $ABDC$, we can suppose WLOG that $E$ is adjacent to the segment $CD$. But this defines the point set up to isometry, and the other distances are $\sqrt{3}$ (two times) and $2$ (one time), so there are no other equilateral triangles among them (the $1$ distances do not contain another either).

Now let us suppose that there are at least 3 of them, which have side length $d$.

For this, we can take the inverse point set, and it leads back to the previous case.

If there are two triangles from both kind, then there are at least $5$ $1$-edges and at least $5$ $d$-edges, so the only way to achieve this would be that there are exactly $5$ edges of both types and our $(1,d)$-graph is the complete graph. But then the $1$-triangles have an edge in common, so part of the point set looks like in Figure \ref{4csucs5el}, but this has only $4$ points and a maximum degree $3$. So the $d$-graph should be isomorphic to the $1$-graph, since that is the only possibility for that too, but on the other hand, it should have a vertex with degree $4$, since the fifth vertex has all $d$-edges. But this is impossible, so we have proved the statement of the lemma for all cases.
\end{pf}

And from now it follows that out of the at least $4$ triangles in the $(1,d)$-graph at least one has sides $1$, $1$, $d$ or $1$, $d$, $d$. And using Proposition \thmref{1dsokszog} for $n=2$, in the former case, the graph only can exist for $d\le 2$, while in the latter case, the graph only can exist for $d\ge{1\over2}$.
\end{pf}

\begin{cor}
Among $5$ points in the plane, there cannot be $4$ different equilateral triangles.
\end{cor}

\begin{pf}
If they have at most $2$ different side lengths, then it immediately follows from Lemma \thmref{szabalyos3szogek}, since there would be such a point set, we could scale it so that one of the side lengths of the equilateral triangles would be $1$. The other case follows from this lemma:

\begin{lem}
Among $5$ points in the plane, equilateral triangles with more than $2$ side lengths cannot exist.
\end{lem}

\begin{pf}
Let us choose one representative of three different side lengths, and all three pairs from these three triangles have at most one vertex in common, otherwise they would have a common side, which is impossible because of their different side lengths. But this means that they have at least $3\cdot3-3\cdot1=6$ vertices in total (without multiplicity), which is also impossible.
\end{pf}
\end{pf}

\begin{prop}
There are no $(1,d)$-graphs with $6$ vertices and $13$ edges for a general $d$.
\end{prop}

\begin{pf}
Such a graph would only have $2$ pairs of vertices which do not have an edge between them, so if we choose one of the endpoints of the first one and a different endpoint of the second one and only take the subgraph induced by the remaining vertices, we get a complete graph with $4$ vertices that is a $(1,d)$-graph. Let us call these $4$ vertices $A$, $B$, $C$ and $D$, while the other two vertices $E$ and $F$. There are $3$ pairs of disjoint edges in this graph. Let such a pair be called a $(1,1)$-type pair, if it contains two $1$-edges, a $(1,d)$-type pair, if it contains a $1$-edge and a $d$-edge and a $(d,d)$-type pair, if it contains two $d$-edges.

\begin{center}
\begin{minipage}{.16\textwidth}
	\centering
	\input{Figures/1dgraf1}
	\captionof{figure}{}
	\label{1dgraf1}
\end{minipage}%
\begin{minipage}{.16\textwidth}
	\centering
	\input{Figures/1dgraf2}
	\captionof{figure}{}
	\label{1dgraf2}
\end{minipage}
\begin{minipage}{.16\textwidth}
	\centering
	\input{Figures/1dgraf3}
	\captionof{figure}{}
	\label{1dgraf3}
\end{minipage}
\begin{minipage}{.16\textwidth}
	\centering
	\input{Figures/1dgraf4}
	\captionof{figure}{}
	\label{1dgraf4}
\end{minipage}
\begin{minipage}{.16\textwidth}
	\centering
	\input{Figures/1dgraf5}
	\captionof{figure}{}
	\label{1dgraf5}
\end{minipage}
\begin{minipage}{.16\textwidth}
	\centering
	\input{Figures/1dgraf6}
	\captionof{figure}{}
	\label{1dgraf6}
\end{minipage}
\end{center}

Red edges mean $1$-edges and blue edges mean $d$-edges.

There are the following two cases:

Among the $3$ pairs, there is at most one $(1,d)$-type pair. This means that there are at least two pairs of the other type, so if we choose those, we either get a parallelogram or an antiparallelogram. If it is a parallelogram with two $(1,1)$-type pairs as its sides, then it has at least one $d$-edge as a diagonal, since its diagonal between its non-obtuse angles is longer than both of the side-lengths. So the other diagonal is either a $1$-edge or a $d$-edge. In the first case we get a rhombus like in Figure \ref{1dgraf1} and in the second case a square like in Figure \ref{1dgraf2}. If it is a parallelogram with two $(d,d)$-type pairs, then its the inverse graph of one of the $(1,d)$-graphs discussed above. If it is a parallelogram with a $(1,1)$-type pair and a $(d,d)$-type pair, then again, the diagonal between the non-obtuse angles should be longer than all of the side-lengths, so it cannot be $1$, neither can it be $d$, which is impossible, so there is no such $(1,d)$-graph. If it is an antiparallelogram, then the only possibility is that it has a $(1,1)$-type pair and a $(d,d)$-type pair, otherwise two of its vertices would be the same. Let us now suppose that $d$ is longer than $1$. The two remaining edges must be of different length, otherwise we would get one of the cases discussed above, since the convex hull of the points would be a parallelogram. So we have got a trapezoid with side lengths $d$, $1$, $1$ and $1$ and diagonal lengths $d$ and $d$. Let us name the four vertices of the trapezoid as it is seen in Figure \ref{1dgraf3}. If we denote $\angle{CBA}=\angle{ACB}$ by $\alpha$ and $\angle{CAD}=\angle{DCA}$ by $\beta$, then $\angle{BAC}=180^\circ-2\alpha$ and $\angle{ADC}=180^\circ-2\beta$, so $180^\circ=\angle{BAD}+\angle{ADC}=\angle{BAC}+\angle{CAD}+\angle{CAD}=180^\circ-2\alpha+\beta+180^\circ-2\beta=360^\circ-2\alpha-\beta$, so $2\alpha+\beta=180^\circ$ and also $\alpha=\angle{BAD}=\angle{BAC}+\angle{CAD}=180^\circ-2\alpha+\beta\Rightarrow3\alpha-\beta=180^\circ$. So $\beta=3\alpha-180^\circ$ meaning that $180^\circ=2\alpha+\beta=5\alpha-180^\circ\Rightarrow\alpha=72^\circ\Rightarrow\beta=3\alpha-180^\circ=36^\circ$. From here, it can be seen that trapezoid $ABCD$ is defined up to isometry and its vertices are $4$ vertices of a regular pentagon.

Among the $3$ pairs, there are at least $2$ $(1,d)$-type pairs. Two of these form a deltoid with side lengths $1$ and $d$ (let $A$ be the point with the two $1$-edges and $C$ be the point with the two $d$-edges). Then either the $AC$ axis of the deltoid has length $1$, and so does $BD$ in which case the only possibility is that the diagonal with length $1$ and the two $1$-sides form an equilateral triangle, while the fourth vertex is at one of the two places shown in \ref{1dgraf4} and \ref{1dgraf5}, or the axis has length $1$, but the other diagonal has length $d$, in which case the non-axis diagonal and the two $d$-sides or the deltoid form an equilateral triangle with side length $d$ and the fourth vertex is its center (\ref{1dgraf6}). There are also the possibilities that both diagonals have length $d$ or that the axis has length $d$ and the other diagonal has length $1$, but these are the inverse graphs of the previously mentioned ones.

So we have finitely many possibilities for the set $\lbrace A,B,C,D\rbrace$ and it also determines $d$, so none of the sets we get from this can have a general $d$ (also, there are only finitely many point sets that can be the representations of such $(1,d)$-graphs since both $E$and $F$ are adjacent to at least two points from the set $\lbrace A,B,C,D\rbrace$, so they are among the intersection points of the circles around these points with radii $1$ and $d$ (no two circles can coincide because any two of them differ in either their radii or their centers or both)).
\end{pf}

\chapter{The chromatic number of the plane}

\section{Introduction}

\begin{defn}
Let us define that graph on the Euclidean plane, whose vertices are the points of the plane and its edges are between the points with distance 1. The chromatic number of this graph will be the \nt{chromatic number of the plane}.
\end{defn}

\begin{thm}
The chromatic number of the plane is equal to the supremum of the chromatic numbers of its finite subgraphs.
\end{thm}

\begin{pf}
The statement trivially follows from the DeBruijn--Erdős theorem.
\end{pf}

\section{Upper bound}

\begin{thm}[description={Isbell, 1950
}]
The chromatic number of the plane is at most 7.
\end{thm}

\begin{center}
	\input{Figures/Isbell}
	\captionof{figure}{}
	\label{isbell}
\end{center}

\begin{pf}
If we periodically extend the colouring on Figure \ref{isbell} to the plane and we choose the side length of the hexagons between ${1\over{\sqrt{7}}}$ and $1\over2$, then we get a good 7-colouring: the diameter of the hexagons will be at most $2\cdot {1\over2}=1$ and the distance of the closest hexagons of the same colour will be at least $\sqrt{\left(5\over2\right)^2+\left({\sqrt{3}\over{2}}\right)^2}\cdot{1\over\sqrt{7}}=\sqrt{{25\over4}+{3\over4}}\cdot{1\over{\sqrt{7}}}=\sqrt{7}\cdot{1\over{\sqrt{7}}}=1$, so there will be no such point pair which has a distance of exactly $1$ and whose members have the same colour. So only the colouring of the borders can cause a problem and if the side lengths are strictly between ${1\over{\sqrt{7}}}$ and ${1\over2}$, then the colour of the borders does not matter. But if for all hexagons we colour the two uppermost vertices and the interior of the three uppermost sides with the same colour as the hexagon, then we get a colouring (for all vertices, there is exactly one hexagon, which has it among its two uppermost vertices and for all sides, there is exactly one hexagon, which has it among its three uppermost sides) and also, this colouring is also good for side lengths ${1\over2}$ and ${1\over{\sqrt{7}}}$ as neither the diameter of a hexagon, nor the shortest distance between two hexagons are included this way.
\end{pf}

\section{Lower bound}

\begin{thm}[description={Nelson, 1950}]
The chromatic number of the plane is at least 4.
\end{thm}

{\rm{\cite{h}}} The following proof is by Hugo Hadwiger:

\begin{pf}
In a good colouring of the plane, there exist $2$ points in the plane, which have distance $\sqrt{3}$, and do not have the same colour, otherwise the whole plane could be easily coloured with one colour as the graph of $\sqrt{3}$ distances is connected, and thus, we would get a contradiction as there are edges in the unit distance graph of the plane. But then we could draw the graph in Figure \ref{4csucs5el} with the two aforementioned points being the acute angles of the rhombus, and we would have to colour the other two vertices with different colours from each other and these two, meaning that we have used four colours for these four points alone, so for the plane too.
\end{pf}

What makes Hadwiger's proof interesting that we can find traces of the idea of the probabilistic method used in Chapter 4 in it. We present here a second proof for the same theorem by William and Leo Moser.

\begin{pf}
The Moser spindle (Figure 1.6) is not 3-colourable (if we try to colour it with $3$ colours ($1$, $2$ and $3$), and $A$ would get colour $1$, then $B$ and $C$ would be coloured with $2$ and $3$ with some order, since the triangle $ABC$ must have all different vertices. The same applies to $D$ and $E$. So $F$ and $G$ both must get colour $1$, which is a contradiction, since they are connected to each other.
\end{pf}

\begin{thm}[description={de Grey, 2018; Exoo, Ismailescu, 2018}]
The chromatic number of the plane is at least 5.
\end{thm}

\begin{pfsketch}
It can be shown by examining the colourings of a graph made from triangular grids that there is an equilateral triangle in the plane with side length $\sqrt{3}$, whose vertices all have the same colour. And from here, he found a graph with $1345$ vertices that cannot have a colouring in which a certain equilateral triangle with side length $\sqrt{3}$ cannot have all of its vertices having the same colour. But the latter part was only verified by computer.
\end{pfsketch}

In the above mentioned result, de Grey have found a not $4$-colourable graph with $20425$ vertices. The Polymath 16 project, which aimed to simplify the prrof resulted in the reduction of the graph. The currently known best result is the following:

\begin{thm}[description={Heule, 2019}]
{\rm{\cite{h}}} There exists a graph with $517$ vertices that is not $4$-colourable.
\end{thm}

\begin{defn}
Let us take a colouring in which the plane is built up from monochromatic regions, which have a positive, existing area and are bordered by Jordan curves. Let us call such a colouring a \nt{tile-based colouring} or a \nt{tiling} of the plane.
\end{defn}

\begin{thm}[description={Townsend, 2005}]
{\rm{\cite{t}}} Every tiling of the plane has at least $6$ colours.
\end{thm}

\chapter{The probabilistic formulation of the Hadwiger--Nelson problem}

\section{Introduction}

The results in this chapter are partly from ${\rm{\cite{pwp}}}$, ${\rm{\cite{pb1}}}$ and ${\rm{\cite{pb2}}}$, but the constructions using the method in Proposition \thmref{5csucsukorlatok} are from my joint work with my supervisor.

We can formulate the question in a probabilistic way:

Let us suppose that we have found a $4$-colouring of the plane such that no two points with distance $1$ have the same colour. Let us create new colourings by composing this colouring with a permutation $\sigma\in S_4$ of the colours on the left and a Euclidean isometry $T\in E(2)$ on the right. This way we created a new colouring $\sigma\circ c\circ T^{-1}:\mathbb C\rightarrow\lbrace 1,2,3,4\rbrace$ of the complex plane, which also has the property of not having a monochromatic pair of points with distance $1$.

The group $S_4\times E(2)$ is solvable, and it is known that solvable groups are amenable, so it is also amenable. So there is a finitely additive probability measure $\mu$ on $S_4\times E(2)$, where all of the subsets of this group are measureable, which is left-invariant: $\mu(gE)=\mu(E)$ for all $g\in S_4\times E(2)$	 and $E\subset S_4\times E(2)$. So $S_4\times E(2)$ has the structure of a finitely additive probability space. Now we can define a random colouring $c:\mathbb C\rightarrow\lbrace 1,2,3,4\rbrace$ by defining $c:=\sigma\circ c\circ T^{-1}$, where $(\sigma,T)$ is an element of $S_4\times E(2)$. For every complex number $z$, the random colour $c(z)$ is a random variable, which takes its values from $\lbrace1,2,3,4\rbrace$. The measure being left-invariant implies that for any $(\sigma,T)\in S_4\times E(2)$, the colouring $\sigma\circ c\circ T^{-1}$ satisfies the condition of not having monochromatic unit distance pairs. From this, we get the colour permutation invariance $\mathbb P(c(z_1)=c_1,...,c(z_k)=c_k)=\mathbb P(c(z_1)=\sigma(z_k)=\sigma(c_k))$ for any $z_1,...,z_k\in\mathbb C, c_1,...,c_k\in\lbrace 1,2,3,4\rbrace$, and $\sigma\in S_4$, and the Euclidean isometry invariance $\mathbb P( c(z_1)=c_1,...,c(z_k)=c_k)=\mathbb P(c(T(z_1))=c_1,...,c(T(z_k))=c_k$.

\begin{defn}
Let us call the event, when points $A_1,A_2,...,A_n$ have the same colour (or in other words, $A_1A_2...A_n$ is monochromatic) $M(A_1A_2...A_n)$.
\end{defn}

\begin{defn}
For some number $d\ge 0$, denote the probability of two points with distance $d$ being monochromatic in the above defined probability field by $p_d$. So $P(M(AB))=p_d$ for any two points $A$ and $B$ with distance $d$.
\end{defn}

This is not necessarily unique, but we can give upper and lower bounds for it. The goal is to find a contradiction by giving a smaller upper bound than the lower bound for some $d$.

\section{Upper bounds}

\begin{prop}[label=d12pd12]
If $d\ge{1\over2}$, then $p_d\le{1\over2}$.
\end{prop}

\begin{pf}
Let us take an isosceles triangle with base-length $1$ and leg-length $d$. In any good $4$-colouring this triangle has at most one monochromatic leg, since if it would have two, its base would also be monochromatic, which is a contradiction. So indeed, $p_d\le{1\over2}$.
\end{pf}

\begin{prop}[label=prop361]
{\rm{\cite{pwp}}} If $d\ge{2\over\sqrt{15}}=0.5163...$, then $p_d\le{1\over2}-{1\over62}$.
\end{prop}

\begin{center}
\begin{minipage}{.49\textwidth}
	\centering
	\input{Figures/Prop361a}
	\captionof{figure}{}
	\label{prop361a}
\end{minipage}
\begin{minipage}{.49\textwidth}
	\centering
	\input{Figures/Prop361b}
	\captionof{figure}{}
	\label{prop361b}
\end{minipage}%
\end{center}

\begin{pf}
Let us take an isosceles triangle $ABC$ with legs $AB$ and $AC$ with length $d$ and with base $BC$ with length $1$. Let us reflect $B$ to $AC$, and name the reflection $D$. Let us now reflect $B$ to the midpoint of $AD$ and let us name the reflection $E$, and let us reflect $A$ to the midpoint of $DE$ and let us name the reflection $F$. Such point sets are drawn in Figures \ref{prop361a} and \ref{prop361b}.

\begin{lem}
If the above point set is drawn with $d\ge{2\over\sqrt{15}}$, then $\abs{AE}$, $\abs{BD}$, $\abs{BE}$, $\abs{BF}$ and $\abs{DF}$ all are at least ${1\over2}$.
\end{lem}

\begin{pf}
Since $B$ and $D$ are each others reflections to $AC$, $ABCD$ is a deltoid. So its area is the sum of the areas of triangles $ABC$ and $ACD$, which is ${\abs{AB}\cdot\abs{BC}\cdot\sin{\angle{CBA}}\over2}+{\abs{AC}\cdot\abs{CD}\cdot\sin{\angle{DCA}}\over2}={d\cdot1\cdot{\sqrt{d^2-{1\over4}}\over d}\over2}+{d\cdot1\cdot{\sqrt{d^2-{1\over4}}\over d}\over2}=\sqrt{d^2-{1\over4}}$. But also, it is 

If $\abs{BD}=x$, then $\sin{\angle{DCB}\over2}={x\over2}$. And since ${\angle{DCB}\over2}={\angle{ACB}}$, which is a base-angle of an isosceles triangle, it is smaller than $90^\circ$, so $\cos{\angle{DCB}\over2}=\sqrt{1-\left({x\over2}\right)^2}$, so $\sin{\angle{DCB}}=2\cdot\sin{\angle{DCB}\over2}\cdot\cos{\angle{DCB}\over2}=2\cdot\frac{x}{2}\cdot\sqrt{1-\left(\frac{x}{2}\right)^2}=\sqrt{x^2-\frac{x^4}{4}}$, so from the sine law we get that $d={x\over{2\cdot\sqrt{x^2-\frac{x^4}{4}}}}$, which is ${{1\over2}\over{2\cdot\sqrt{\frac{1}{4}-\frac{1}{64}}}}=\frac{2}{\sqrt{15}}$ for $x=\frac{1}{2}$ and $d$ is strictly monotonously increasing by $x$, since $d={x\over{2\cdot\sqrt{x^2-\frac{x^4}{4}}}}={1\over{2\cdot\sqrt{1-{x^2}\over4}}}$, so for $x<{1\over2}$, $d<{2\over\sqrt{15}}$, which confirms that from $d\ge{2\over\sqrt{15}}$, $\abs{BD}\ge{1\over2}$ follows. And because of the isometry of triangles $ABD$ and $DEA$, $\abs{EA}\ge{1\over2}$ also follows and because of the isometry of triangles $ADE$ and $FED$, $\abs{DF}\ge{1\over2}$ is also true. And because of triangle $ABD$ being isosceles, $\angle{DBA}$ is an acute angle (although its orientation is ambigous), so in parallelogram $ABDE$, $BE$ is a longer diagonal, so as sides all sides have length at least ${1\over2}$, it is at least ${1\over2}$ too. Similarly $\angle{DAE}$ is an acute angle too, so parallelogram $ADFE$ has $AF$ as a longer diagonal, so $\abs{AF}\ge{1\over2}$ as the sides of this parallelogram also have lengths at least ${1\over2}$. Also, since $\vec{BD}=\vec{AE}$ because of the isometry of triangles $ABD$ and $DEA$, and $\vec{AE}=\vec{DF}$ because of the isometry of triangles $DEA$ and $EDF$, we get $\vec{BF}=\vec{BD}+\vec{DF}=2\cdot\vec{BD}$, which means that it also has length at least $1>{1\over2}$.
\end{pf}

Let us now define $\delta$ as ${1\over2}-p_d$. This is a non-negative number because of \thmref{d12pd12}. Now $AB$ and $AC$ both have a probability $P(M(AB))=P(M(AC))=p_d={1\over2}-\delta$ for being monochromatic, and since both cannot be monochromatic at the same time as that would mean $BC$ being monochromatic, $M(AB)$ and $M(AC)$ are disjoint events, so $P(\overline{M(AB)}\cap\overline{M(AC)})=2\delta$. The same argument applies to $M(AC)$ and $M(AD)$. And since $AB$ and $AC$ cannot be both monochromatic (as it was stated above), $P(M(AD))\subseteq P(\overline{M(AC)})\Rightarrow P(M(AD)\cap\overline{M(AB)})\le P(\overline{M(AC)}\cap\overline{M(AB)})=2\delta$. The same argument can be applied to the other two triangles isometric to $ABD$: $P(M(DE)\cap\overline{M(DA)})\le2\delta$ and $P(M(EF)\cap\overline{M(ED)})\le2\delta$.

Let us now define the event $\epsilon=M(AB)\cup M(AD)\cup M(DE)\cup M(EF)$. From the above, $P(M(AB))=\frac{1}{2}-\delta$, $P(M(AD)\setminus M(AB))\le2\delta$, $P(M(DE)\setminus M(AD))\le2\delta$ and $P(M(EF)\setminus M(DE))\le2\delta$, so as these events together cover $\varepsilon$, $P(\varepsilon)\le{1\over2}-\delta+2\delta+2\delta+2\delta={1\over2}+5\delta$. And since $P(M(AD))={1\over2}-\delta$ and $P(M(AD)\setminus M(AB))\le2\delta$, $P(M(ABD))=P(M(AB)\cap M(AD))=P(M(AD))-P(M(AD)\setminus M(AB))\ge\frac{1}{2}-3\delta$. But $M(ABD)$ is in $\varepsilon$ and, because of the length of $BD$, $P(M(BD))\le\frac{1}{2}$, thus, $P(M(BD)\setminus\varepsilon)$ is at most $P(M(BD))-P(M(ABD))\le3\delta$. Similarly $P(M(AE)\setminus\varepsilon)$ is also at most $3\delta$ as is $P(M(DF)\setminus\varepsilon)$.

Also, $P(M(ABDE))\ge P(M(DE))-P(M(DE)\setminus M(AD))-P(M(AD)\setminus M(AB))\ge\frac{1}{2}-5\delta$, and $M(ABDE)\subseteq\varepsilon$ so because of $P(M(EB))\le\frac{1}{2}$, $P(M(EB)\setminus\varepsilon)\le P(M(EB)-P(M(ABDE))\le\frac{1}{2}-5\delta$. A similar argument applies to $M(AF)$.

And finally $P(M(BF)\setminus\varepsilon)\le7\delta$ from a similar argument, so the probability that at least one pair from $A,B,D,E,F$ is monochromatic is at most $\frac{1}{2}+5\delta+3\delta+3\delta+3\delta+5\delta+5\delta+7\delta={1\over2}+31\delta$. On the other hand, from the pigeonhole principle, this probability is $1$, proving $31\delta\ge\frac{1}{2}$ and thus finishing the statement.
\end{pf}

\begin{prop}
If $d\ge{\sqrt{3}-1\over\sqrt{2}}=0.5176...$, then $p_d\le0.48$.
\end{prop}

\begin{pf}
Let us take the same point set as in \thmref{prop361}. We now not only want $\abs{BD}$ to be at least ${1\over2}$, but also to be at least $d$. $\abs{BD}=d$ if and only if $ABD$ is an equilateral tirangle, which happens if and only if the isosceles triangles $ABC$ and $ACD$ have apex angle $150^\circ$ or $30^\circ$. The first one occurs if $d={1\over2\sin{75^\circ}}={\sqrt{3}-1\over\sqrt{2}}\approx 0.5176$, while the second one occurs if $d={1\over2\sin{15^\circ}}={\sqrt{3}+1\over\sqrt{2}}\approx1.9318$. Between these values $\abs{BD}>d$, since it is a chord of the circle around $A$ with radius $d$ that belongs to a larger central angle than $60^\circ$ (and smaller than $300^\circ$). And this means that $d\ge{1\over2\sin(75^\circ)}={\sqrt{3}-1\over\sqrt{2}}=0.5176...$ and $d\le{1\over2\sin(15^\circ)}={\sqrt{3}+1\over\sqrt{2}}=1.9318...$, but the latter is not needed, it is enough that $\abs{BD}\ge d_0={\sqrt{3}-1\over\sqrt{2}}$, which is true for all such cases. If we pick a $d\ge d_0$ for which $p_d\ge{1\over2}-\delta-\varepsilon$, where $\sup_{d\ge d_0}{p_d}={1\over2}-\delta$ and $\varepsilon$ is a small positive number. If we do the same calculation as in the first case, we get ${1\over2}+5\delta+2\delta+2\delta+2\delta+4\delta+4\delta+6\delta+O(\varepsilon)={1\over2}+25\delta+O(\varepsilon)\ge 1$ meaning that $\delta\ge 0.02$ if we have chosen an $\varepsilon$ small enough.
\end{pf}

\begin{prop}
$\limsup{d}{p_d}\le{323\over675}=0.4785...$ (in other words, $p_d<0.4786$ if $d$ is large enough).
\end{prop}

\begin{pfsketch}
The proof is similar to the previous ones, but uses a slightly different point set: now we define $F$ by reflecting $B$ to $AD$ and we tend to infinity with $d$.
\end{pfsketch}

\begin{prop}
If $d>{1\over2}$, $p_d<{1\over2}$.
\end{prop}

\begin{pf}
Let us take a circle of radius $d$ centered at $O$ and an infinite series of points $A_1$, $B_1$, $A_2$, $B_2$, $A_3$, $B_3$, $A_4$, $B_4$, ... on the circle so that the distance of the neighbouring points is $1$ (we can do this with the exception of those numbers, which are the side lengths of regular (possibly star) polygons on the unit circle, which case we will discuss later). Then assuming $p_d=0.5$ (it cannot be larger because of \thmref{d12pd12}), then either all the $A_i$'s or the $B_i$'s (or at least the first $n$ of them for an arbitrarily large $n$) have the same colour as $O$ with probability $1$ since all the radii ending in the $A_i$'s and the $B_i$'s have length $d$ and the neighbouring points in the $A_1B_1A_2B_2...$ sequence cannot have the same colour. We can choose four of the $A_i$'s which all have distance larger than $d$, since distance $d$ belongs to arcs with central angle $60^\circ$ and we can find an arbitrarily small arc among the $A_i$'s and we also can step with it as many as we want. Let us call these four $A_i$'s $A_{i_1}$, $A_{i_2}$, $A_{i_3}$ and $A_{i_4}$. As it was stated above, either the $A_i$'s or the $B_i$'s are monochromatic, so it has probability at least ${1\over2}$, that all the $A_i$'s are monochromatic, meaning that the $A_{i_j}$'s (for $j\in\lbrace1,2,3,4\rbrace$) are also monochromatic along with $O$. But it means that if the $B_i$'s are the monochromatic ones, than the points $A_{i_1}$, $A_{i_2}$, $A_{i_3}$, $A_{i_4}$ and $O$ have different colours (since for all distances $\ge{1\over2}$, $p_d\le{1\over2}$), which is impossible, because there are only $4$ colours. So we have got a contradiction, thus the statement follows.

For the regular polygon case: $\arccos{\left({1\over2}\right)}=60^\circ$, so for enough sides, we can find four points like it was described above, while for the small cases, it is already known that $p_d<{1\over2}$, since even a star polygon (with few sides) cannot have such a small angle, which would make $d$ to be closer to ${1\over2}$, than $\approx0.5176$.
\end{pf}

\begin{rem}
A similar argument can be made if we do not use a circle, but instead we use a chain of isosceles triangles. The advantage of using such a method is that it can be used for arbitrarily many colours.
\end{rem}

\begin{rem}
With another similar argument, we can also see that if $x$, $y$ and $1$ are the sides of a non-degenerate triangle, then $p_x+p_y<1$.
\end{rem}

\section{Lower bounds}

\begin{thm}[label=5pont]
Let us take 5 different points such that among the distances between them, there are at least 8, which are $1$ or $d$ for some fixed $d$. Let us call the two remaining distances $x$ and $y$ in increasing order (or if only one distance remained, let's call it $x$). Then $p_d\ge{1-p_x-p_y\over n}$ (or in case there is only 1 non-$1$, non-$d$ distance, $p_d\ge{1-p_x\over n}$, and if there is no such distance, $p_d\ge {1\over n}$) (where $n$ denotes the number of distances with length $d$). We are only interested in such proofs if they give a positive lower bound for $p_d$.
\end{thm}

\begin{defn}
If for a $5$ element point set the method described above gives us a proof that $p_d$ is positive, let us call that point set a \nt{good $5$ point set}.
\end{defn}

\begin{pf}
There are at least two points of the same colour among the 5 points, so the expected number of monochromatic pairs is at least 1. But all of the monochromatic pairs should be among those, which have distance $d$, $x$ or $y$ (since those with distance 1 can't be of the same colour). So $n\cdot p_d+p_x+p_y\ge1$ (or in case of 1 non-$1$, non-$d$ distance, $n\cdot p_d+p_x\ge1$ and in case of no such distances, $n\cdot p_d\ge$), from which the statement follows.
\end{pf}

\begin{thm}[label=5csucsukorlatok]
With the known upper bounds, we cannot give an analogous proof for any lower bound for any $d$ if there are less than 8 distances from $\lbrace 1,d\rbrace$.
\end{thm}

\begin{pf}
The minimal known upper bound is ${1\over3}$, so the sum of the three upper bounds for the remaining distances is at least 1, so $1-p_x-p_y-p_z\le 0$ (where $x$, $y$ and $z$ are the three non-$1$, non-$d$ distances). So the lower bound we would get for $p_d$ is at most $0$.
\end{pf}

The following are the bounds, which can be concluded from \thmref{5csucsukorlatok} using \thmref{nehany1dgraf}.

\begin{prop}[label=5csucsugraf1]
$p_{\sqrt{6}+\sqrt{2}\over2}\ge{13\over75}$.
\end{prop}

\begin{pf}
Let us take graph 6 from Example \thmref{nehany1dgraf} $BAE$, $EAC$ and $DAB$ are all isosceles triangles with side length $1$, apex $A$ and a $150^\circ$ angle in $A$, so $BE$, $EC$ and $DB$ are all $\sqrt{1^2+1^2-2\cdot1\cdot1\cdot\left(-{\sqrt{3}\over2}\right)}=\sqrt{2-\sqrt{3}}={\sqrt{6}+\sqrt{2}\over2}$ long. $DAC$ is also an isosceles triangle with apex $A$, but here the angle in $A$ is $90^\circ$, so $CD$ is $\sqrt{2}$. So among the five points there are $6$ segments of length $1$, $3$ segments of length ${\sqrt{6}+\sqrt{2}\over2}$ and $1$ segment of length $\sqrt{2}$ for $d={\sqrt{6}+\sqrt{2}\over2}$, $p_d\ge{1-p_x\over3}={1-p_{\sqrt{2}}\over3}\ge{1-{12\over25}\over3}={13\over75}$.
\end{pf}

\begin{prop}[label=5csucsugraf2]
$p_{\sqrt{6}-\sqrt{2}\over2}\ge{13\over150}$.
\end{prop}

\begin{pf}
Let us take the inverse graph of graph 6 from Example \thmref{nehany1dgraf} and from \thmref{5csucsugraf1} we can see that $d={\sqrt{6}-\sqrt{2}\over2}$ and there are $6$ $d$-edges and $3$ $1$-edges, which proves the statement from \thmref{5csucsukorlatok}.
\end{pf}

\begin{prop}[label=5csucsugraf3]

$d\le{\sqrt{6}-\sqrt{2}\over2}\Rightarrow p_d\ge{1\over50}$

${\sqrt{6}-\sqrt{2}\over2}<d\le\sqrt{14\over15}-{1\over\sqrt{15}}\Rightarrow p_d\ge{14\over775}$

$\sqrt{14\over15}-{1\over\sqrt{15}}<d\le{\sqrt{15}\over4}-{\sqrt{3}\over4}\Rightarrow p_d\ge{1\over100}$

${\sqrt{15}+\sqrt{3}\over4}\le d<\sqrt{14\over15}+{1\over\sqrt{5}}\Rightarrow p_d\ge{1\over100}$

$\sqrt{14\over15}+{1\over\sqrt{5}}\le d<\sqrt{2}\Rightarrow p_d\ge{14\over775}$

$\sqrt{2}\le d<2\Rightarrow p_d\ge{1\over50}$
\end{prop}

\begin{center}
	\input{Figures/5pontu1dgraf1}
	\captionof{figure}{}
	\label{moser}
\end{center}

\begin{pf}
Let us take graph number 1 from Example \thmref{nehany1dgraf} and let us suppose that $\varphi\in[0,{\pi\over2}]$.

$\abs{CD}\le\abs{BE}$, since they are both chords of the unit circle around $A$, but $CD$ belongs to a central angle of size $\varphi-60^\circ$, while $BE$ belongs to an angle of size $300^\circ-\varphi$ and since the difference of the angles is $(300^\circ-\varphi)-(\varphi-60^\circ)=360^\circ-2\cdot\varphi$ and $\varphi\le 180^\circ$, $\angle EAB\ge\angle CAD$. This means that $x=\abs{CD}$ and $y=\abs{BE}$.

From the above formula it is also clear that $\angle EAB$ is at least $120^\circ$ and at most $300^\circ$, so $\cos{\angle EAB}\ge{1\over2}$, thus from the cosine theorem, $\abs{BE}^2\ge1^2+1^2-2\cdot1\cdot1\cdot\left({1\over2}\right)=1$, which means that $y=\abs{BE}\ge 1\Rightarrow p_y\ge{12\over25}$.

Let us now take a coordinate system, where the origin is the midpoint of $C$ and $D$, the coordinates of $C$ are $\left(-{x\over2};0\right)$, the coordinates of $D$ are $\left({x\over2};0\right)$ (and from now on, if triangles $ABC$ and $ADE$ cross each other, we define $x$ as $-\abs{CD}$).

We can see that $x$ moves between $-1$ and $\sqrt{3}$, because its value is $2\cdot\sin{\left({\varphi-60^\circ\over2}\right)}$, where ${\varphi-60^\circ\over2}$ moves in the interval $\left[-30^\circ,60^\circ\right]$. But similarly, $y$ moves between $\sqrt{3}$ and $1$, so it is always larger than ${\sqrt{3}-1\over\sqrt{2}}$.

The coordinates of $A$ are $\left(0,\sqrt{1-{x^2\over4}}\right)$, so the coordinates of the midpoint of $AD$ (let us call it $M$) are $\left({x\over4},{\sqrt{1-{x^2\over4}}\over2}\right)$. And since $\vec{ME}$ is the rotated image of $\vec{MD}$ with $90^\circ$ scaled with $\sqrt{3}$, $E=\left({x\over4},{\sqrt{1-{x^2\over4}}\over2}\right)+\left({\sqrt{3-{3x^2\over4}}\over2},{\sqrt{3}\cdot x\over4}\right)=\left({x\over4}+{\sqrt{3-{3x^2\over4}}\over2},{\sqrt{1-{x^2\over4}}\over2}+{\sqrt{3}\cdot x\over4}\right)$, so
\begin{align*}
	d & =\abs{CE}
		=\sqrt{\left({3\over4}x+{\sqrt{3-{3x^2\over4}}\over2}\right)^2
		+\left({\sqrt{1-{x^2\over4}}\over2}+{\sqrt{3}\cdot x\over4}\right)^2}= \\
	& =\sqrt{{9\over16}x^2+{12-3x^2\over16}+{3\over4}\cdot\sqrt{3-{3x^2\over4}}\cdot x+{4-x^2\over16}+{3x^2\over16}+{1\over4}\cdot\sqrt{3-{3x^2\over4}}}= \\
	& =\sqrt{{x^2\over2}+1+\sqrt{3-{3x^2\over4}\cdot x}}
\end{align*}
(and the range of this function indeed includes $\left(-1,\sqrt{3}\right]$), so we get the statement of the proposition.
\end{pf}

\begin{prop}[label=5csucsugraf4]

${1\over2}\le d\le\sqrt{{1\over2}\left(2-\sqrt{2}\right)\left(3-\sqrt{3}\right)}\Rightarrow p_d\ge{1\over150}$

$\sqrt{{1\over2}\left(2-\sqrt{2}\right)\left(3-\sqrt{3}\right)}<d\le\sqrt{{19\over15}-{2\over\sqrt{15}}}\Rightarrow p_d\ge{14\over2325}$

$\sqrt{{19\over15}-{2\over\sqrt{15}}}<d\le{1\over2}\sqrt{5-2\sqrt{3}}\Rightarrow p_d\ge{1\over300}$

${1\over2}\sqrt{5+2\sqrt{3}}\le d<\sqrt{{19\over15}+{2\over\sqrt{15}}}\Rightarrow p_d\ge{1\over300}$

$\sqrt{{19\over15}+{2\over\sqrt{15}}}\le d<\sqrt{{1\over2}\left(2-\sqrt{2}\right)\left(3-\sqrt{3}\right)}\Rightarrow p_d\ge{14\over2325}$

$\sqrt{{1\over2}\left(2+\sqrt{2}\right)\left(3-\sqrt{3}\right)}\le d\Rightarrow p_d\ge{1\over150}$
\end{prop}

\begin{center}
	\input{Figures/5pontu1dgraf2}
	\captionof{figure}{}
	\label{moser}
\end{center}

\begin{pf}
Let us take the $(1,d)$-graph number 2 from Example \thmref{nehany1dgraf} and let us suppose (we can do that WLOG) that $\angle{CAD}\le120^\circ$. This way again $\abs{CD}\le\abs{EB}$ (since both are chords of the circle around $A$ with radius $d$, and $EB$ belongs to the larger central angle), so $x=\abs{CD}$. This time $x$ will not be a signed value. Let us draw the quadrilateral $ACDE$ so that $C$ is the origin and $E$ is the $(1,0)$ point. This way $A$ is the apex of an isosceles triangle with base $CE$ and legs of length $d$, so its coordinates are $\left({1\over2},\sqrt{d^2-{1\over4}}\right)$. The coordinates of the midpoint of $A$ and $E$ are $\left({3\over4},{\sqrt{d^2-{1\over4}}\over2}\right)$ so the coordinates of $D$ (which is the third vertex of an equilateral triangle with side $AE$) are $\left({3\over4}-{\sqrt{3d^2-{3\over4}}\over2},{\sqrt{d^2-{1\over4}}\over2}-{\sqrt{3}\over4}\right)$, so the $x=\abs{CD}=\sqrt{\left({3\over4}-{\sqrt{3d^2-{ 3\over4}}\over2}\right)^2+\left({\sqrt{d^2-{1\over4}}\over2}-{\sqrt{3}\over4}\right)^2}=\sqrt{4\cdot\left({\sqrt{d^2-{1\over4}}\over2}-{\sqrt{3}\over4}\right)^2}=2\cdot\left\lvert{\sqrt{d^2-{1\over4}}\over2}-{\sqrt{3}\over4}\right\rvert=\left\lvert\sqrt{d^2-{1\over4}-{\sqrt{3}\over2}}\right\rvert$. So $d=\sqrt{\left(x+{\sqrt{3}\over2}\right)^2+{1\over4}}$ if the triangles intersect each other and $d=\sqrt{\left(x-{\sqrt{3}\over2}\right)^2+{1\over4}}$ if they do not. So in the not crossing case

$$x={\sqrt{3}-1\over \sqrt{2}}\Rightarrow d=\sqrt{{1\over2}\left(2-\sqrt{2}\right)\left(3-\sqrt{3}\right)}\approx 0.609404$$

$$x={2\over\sqrt{15}}\Rightarrow d=\sqrt{{19\over15}-{2\over\sqrt{15}}}\approx 0.610114$$

$$x={1\over2}\Rightarrow d={1\over2}\sqrt{5-2\sqrt{3}}\approx 0.619657$$, while in the crossing case

$$x={1\over2}\Rightarrow d={1\over2}\sqrt{5+2\sqrt{3}}\approx 1.45466$$

$$x={2\over\sqrt{15}}\Rightarrow d=\sqrt{{19\over15}+{2\over\sqrt{15}}}\approx 1.47007$$

$$x={\sqrt{3}-1\over\sqrt{2}}\Rightarrow d=\sqrt{{1\over2}\left(2+\sqrt{2}\right)\left(3-\sqrt{3}\right)}\approx 1.47123$$

From these, and from \thmref{5pont} we can conclude the statement as $n=6$.
\end{pf}

\begin{prop}[label=5csucsugraf5]

${1\over2}\le d\le\sqrt{\left(3-\sqrt{3}\right)\cdot\left(1-{1\over\sqrt{2}}\right)}\Rightarrow p_d\ge{1\over125}$

$\sqrt{\left(3-\sqrt{3}\right)\cdot\left(1-{1\over\sqrt{2}}\right)}<d\le\sqrt{{19\over15}-{2\over\sqrt{5}}}\Rightarrow p_s\ge{28\over3875}$

$\sqrt{{19\over15}-{2\over\sqrt{5}}}<d\le{1\over2}\sqrt{5-2\sqrt{3}}\Rightarrow p_d\ge{1\over250}$

${1\over2}\sqrt{5+2\sqrt{3}}\le d<{1\over2}\sqrt{19+6\sqrt{5}}\Rightarrow p_d\ge{1\over250}$

${1\over2}\sqrt{19+6\sqrt{5}}\le d<\sqrt{\left(3-\sqrt{3}\right)\left(1+{1\over\sqrt{2}}\right)}\Rightarrow p_d\ge{28\over3875}$

$\sqrt{\left(3-\sqrt{3}\right)\left(1+{1\over\sqrt{2}}\right)}\le d\Rightarrow p_d\ge{1\over125}$
\end{prop}

\begin{center}
	\input{Figures/5pontu1dgraf3}
	\captionof{figure}{}
	\label{moser}
\end{center}

\begin{pf}
Let us take the $(1,d)$-graph number 3 from Example \thmref{nehany1dgraf} and let us suppose WLOG that lines $AD$ and $AE$ are rotated with a positive angle compared to $AB$. This way if we choose $AB$ to be horizontal and $ABC$ to be directed upwards and we denote it by $x$ that by what measure is $D$ higher than $C$ (so the length of $x$ with a sign), then $d=\abs{AD}=\sqrt{\left({1\over2}\right)^2+\left({\sqrt{3}\over2}+x\right)^2}=\sqrt{{1\over4}+{3\over4}+x^2+\sqrt{3}x}=\sqrt{1+x^2+\sqrt{3}x}$.
This is obviously a strictly monotonically increasing function and all of the possible values of $x$ ($x\in\left(-{\sqrt{3}\over2},\infty\right)$) are in its domain. And since its values are

$$x=-{\sqrt{3}-1\over\sqrt{2}}\Rightarrow d=\sqrt{\left(3-\sqrt{3}\right)\left(1-{1\over\sqrt{2}}\right)}\approx 0.609404$$

$$x=-{2\over\sqrt{15}}\Rightarrow d=\sqrt{{19\over15}-{2\over\sqrt{5}}}\approx 0.610114$$

$$x=-{1\over2}\Rightarrow d={1\over2}\sqrt{5-2\sqrt{3}}\approx 0.619657$$

$$x={1\over2}\Rightarrow d={1\over2}\sqrt{5+2\sqrt{3}}\approx 1.45466$$

$$x={2\over\sqrt{15}}\Rightarrow d=\sqrt{{19\over15}+{2\over\sqrt{15}}}\approx 1.47007$$

$$x={\sqrt{3}-1\over\sqrt{2}}\Rightarrow d=\sqrt{\left[3-\sqrt{3}\right)\left(1+{1\over\sqrt{2}}\right)}\approx 1.47123$$

From these and from Proposition \thmref{5pont} we can conclude the statement as $n=5$.
\end{pf}

\begin{prop}[label=5csucsugraf6]

$d\le{\sqrt{6}-\sqrt{2}\over2}\Rightarrow p_d\ge{1\over75}$

${\sqrt{6}-\sqrt{2}\over2}<d\le\sqrt{14\over15}-{1\over\sqrt{15}}\Rightarrow p_d\ge{28\over2325}$

$\sqrt{14\over15}-{1\over\sqrt{15}}<d\le{\sqrt{15}\over4}-\sqrt{3}{1\over4}\Rightarrow p_d\ge{1\over150}$

${1\over2}\sqrt{5+2\sqrt{3}}\le d<\sqrt{{19\over15}+{2\over\sqrt{5}}}\Rightarrow p_d\ge{1\over150}$

$\sqrt{{19\over15}+{2\over\sqrt{5}}}\le d<\sqrt{{1\over2}\left(2+\sqrt{2}\right)\left(3-\sqrt{3}\right)}\Rightarrow p_d\ge{28\over2325}$

$\sqrt{{1\over2}\left(2+\sqrt{2}\right)\left(3-\sqrt{3}\right)}\le d\le2\Rightarrow p_d\ge{1\over75}$
\end{prop}

\begin{center}
	\input{Figures/5pontu1dgraf4}
	\captionof{figure}{}
	\label{moser}
\end{center}

\begin{pf}
Now take the inverse case of the first one, so $(1,d)$-graph number 4 from Example \thmref{nehany1dgraf}. Let us again suppose that $ADE$ is rotated compared to $ABC$ with an angle from $\left(0,{\pi\over2}\right]$. Then $D$ is on the unit circle around $B$ and sure not under line $AB$ as that would mean that $\angle{BAD}\ge{\pi\over2}$. If $D$ is on the shorter arc between $A$ and $C$, then $\angle{ADC}=150^\circ$, since the central angle opposite to it is $60^\circ$. So from the cosine-theorem we get $d^2+x^2+\sqrt{3}dx=1$.

So we have got a quadratic equation and its solution is $d={-\sqrt{3}\pm\sqrt{3x^2-4x^2+4}\over2}={-\sqrt{3}x\pm\sqrt{4-x^2}\over2}$, but because $x$ now cannot be negative (it will be a separate case if the figure is arranged differently), the smaller roots are always negative and thus they cannot be good, so $d={-\sqrt{3}x+\sqrt{4-x^2}\over2}$ and $d$ is strictly monotonically decreasing in $x$. So

$$x={\sqrt{3}-1\over2}\Rightarrow d={\sqrt{6}-\sqrt{2}\over2}\approx 0.517638$$

$$x={2\over\sqrt{15}}\Rightarrow d=\sqrt{14\over15}-{1\over\sqrt{15}}\approx 0.518878$$

$$x={1\over2}\Rightarrow d={\sqrt{15}\over4}-\sqrt{3}{1\over4}\approx 0.535233$$

Similarly if $D$ is on the $CF$ arc, then by the cosine theorem applied to triangle $ADC$, we get $d=\sqrt{1+x^2+\sqrt{x}}$, so

$$x={1\over2}\Rightarrow d={1\over2}\sqrt{5+2\sqrt{3}}\approx 1.45466$$

$$x={2\over\sqrt{15}}\Rightarrow d=\sqrt{{19\over15}+{2\over\sqrt{5}}}\approx 1.47007$$

$$x={\sqrt{3}-1\over\sqrt{2}}\Rightarrow d=\sqrt{{1\over2}\left(2+\sqrt{2}\right)\left(3-\sqrt{3}\right)}\approx 1.47123$$

So from these and \thmref{5pont} we can conclude the statement using that $n=3$.
\end{pf}

\begin{prop}[label=5csucsugraf7]
{\rm{\cite{wa2}}}{\rm{\cite{wa3}}}{\rm{\cite{wa4}}}

$$0<d\le  0.51763...\Rightarrow p_d\ge{1\over100}$$

$$0.51763...<d\le0.51814...\Rightarrow p_d\ge{7\over775}$$

$$0.51814...<d\le0.52561...\Rightarrow p_d\ge{1\over200}$$

$$1.3587...\le d<1.36855...\Rightarrow p_d>0\ge{1\over200}$$

$$1.36855...\le d<1.36929...\Rightarrow p_d>0\ge{7\over775}$$

$$1.36929...\le d\ge 2\Rightarrow p_d>0\ge{1\over100}$$
\end{prop}

\begin{center}
	\input{Figures/5pontu1dgraf5}
	\captionof{figure}{}
	\label{moser}
\end{center}

\begin{pfsketch}
It can be proved with the help of graph number 5 from Example \thmref{nehany1dgraf} using that $\abs{EA}=\sqrt{\left({d^2\over2}-{1\over2}\right)^2+\left(\sqrt{d^2-{d^4\over4}}-\sqrt{d^2-{1\over4}}\right)^2}$.
\end{pfsketch}

\begin{prop}
A good $5$ point set cannot be realizable for numbers below $1\over 2$ and above $2$ at the same time.
\end{prop}

\begin{pf}
Above we have proved that there are at least $8$ segments with length from $\lbrace 1,d \rbrace$ in such a point set. From here, the statement follows from Proposition \thmref{1dgraf58}
\end{pf}

\begin{prop}
If $x\ge{d\over2}$, then $2\cdot(1-p_x)\ge1-p_d$.
\end{prop}

\begin{pf}
Let's take an isosceles triangle with base-length $d$ and leg-length $x$. If the two endpoints of the base have different colours, then the apex has a different colour from at least one of them, so at least one of the legs is bichromatic. And from this the statement follows.
\end{pf}

\begin{prop}
If $0\le d\le {\sqrt{5}-1\over\sqrt{3}}$, then $p_d\ge{1\over325}$.
\end{prop}

\begin{pf}
First we prove that $x\ge {d\over2}$ if we choose $x$ from point set number 2 belonging to $d$:

In point set number 1, for $d\ge {\sqrt{3}+\sqrt{15}\over4}$, $x\ge{1\over2}$ and point set number 2 is the inverse of point set number 1, which proves the statement.

Now let us suppose for contradiction that $p_d<{1\over325}$. From this: $p_d<{1\over325}\Rightarrow 1-p_d>{324\over325}\Rightarrow 1-p_x>{162\over325}\Rightarrow 1-p_x-p_y>{6\over325}\Rightarrow p_d>{1\over325}$, which is a contradiction.
\end{pf}

\begin{prop}
Summarizing the above bounds, we get

$$d\in\left(0,{\sqrt{6}-\sqrt{2}\over2}\right]\Rightarrow p_d\ge{1\over50}$$ (from \thmref{5csucsugraf3})

$$d\in\left({\sqrt{6}-\sqrt{2}\over2}\right]\Rightarrow p_d\ge{14\over775}$$ (from \thmref{5csucsugraf3})

$$d\in\left(\sqrt{14\over15}-{1\over\sqrt{15}},{\sqrt{15}-\sqrt{3}\over4}\right]\Rightarrow p_d\ge{1\over100}$$ (from \thmref{5csucsugraf3})

$$d\in\left({\sqrt{15}-\sqrt{3}\over4},\sqrt{\left(3-\sqrt{3}\right)\left(1-{1\over\sqrt{2}}\right)}\right]\Rightarrow p_d\ge{1\over125}$$ (from \thmref{5csucsugraf5})

$$d\in\left(\sqrt{\left(3-\sqrt{3}\right)\left(1-{1\over\sqrt{2}}\right)},\sqrt{{19\over15}-{2\over\sqrt{5}}}\right]\Rightarrow p_d\ge{38\over3875}$$ (from \thmref{5csucsugraf5})

$$d\in\left(\sqrt{{19\over15}-{2\over\sqrt{5}}},{1\over2}\sqrt{5-2\sqrt{3}}\right]\Rightarrow p_d\ge{1\over250}$$ (from \thmref{5csucsugraf5})

$$d\in\left[1.3587...,1.36855...\right)\Rightarrow p_d\ge{1\over200}$$ (from \thmref{5csucsugraf7})

$$d\in\left[1.36855...,1.36929\right)\Rightarrow p_d\ge{7\over775}$$ (from \thmref{5csucsugraf7})

$$d\in\left[1.36929...,\sqrt{14\over15}+{1\over\sqrt{15}}\right)\Rightarrow p_d\ge{1\over100}$$ (from \thmref{5csucsugraf7})

$$d\in\left[\sqrt{14\over15}+{1\over\sqrt{15}},\sqrt{2}\right)\Rightarrow p_d\ge{14\over775}$$ (from \thmref{5csucsugraf1})

$$d\in\left[\sqrt{2},2\right]\Rightarrow p_d\le{1\over50}$$ (from \thmref{5csucsugraf1})

$$d\in\left(2,\infty\right]\Rightarrow p_d\ge{1\over125}$$ (from \thmref{5csucsugraf5})
\end{prop}

\begin{defn}
Let us take a set of points and let us colour them arbitrarily with $4$ colours and let us call the minimal possible number of monochromatic point pairs $f(n)$.
\end{defn}

\begin{prop}
$f(n)={\left({\left\lfloor{n\over4}\right\rfloor}\choose2\right)}+{\left({\left\lfloor{n+1\over4}\right\rfloor}\choose2\right)}+{\left({\left\lfloor{n+2\over4}\right\rfloor}\choose2\right)}+{\left({\left\lfloor{n+3\over4}\right\rfloor}\choose2\right)}$.
\end{prop}

\begin{pf}
We can easily find such a colouring, where this is the number of monochromatic pairs: we colour ${\left\lfloor{n\over4}\right\rfloor}$ points with one colour, ${\left\lfloor{n+1\over4}\right\rfloor}$ points with another, ${\left\lfloor{n+2\over4}\right\rfloor}$ points with a third colour and ${\left\lfloor{n+3\over4}\right\rfloor}$ points with a fourth one. And as the fractional parts of the numbers ${n\over4}$, ${n+1\over4}$, ${n+2\over4}$ and ${n+3\over4}$ are $0$, ${1\over4}$, ${1\over2}$ and ${3\over4}$ in some order, ${\left\lfloor{n\over4}\right\rfloor}+{\left\lfloor{n+1\over4}\right\rfloor}+{\left\lfloor{n+2\over4}\right\rfloor}+{\left\lfloor{n+3\over4}\right\rfloor}=\left({n\over4}+{n+1\over4}+{n+2\over4}+{n+3\over4}\right)-\left(0+{1\over4}+{1\over2}+{3\over4}\right)={n+6\over4}-{3\over2}=n$. Also, let us suppose that there is some other distribution of points with less monochromatic pairs. But this is the only way that $n$ can be written as the sum of $4$ integers with a maximal difference $1$ (the other option would be that there are other four numbers of form ${k\over4}$, whose integer parts sum up to $n$, but we can choose them to be consecutive and depending on they are larger or smaller than the above four, their sum is larger or smaller too, since in every step, the sum grows with $1$. And if there would be two numbers $a$ and $b$ with $a-b>1$, then we could reduce the number of monochromatic pairs: ${\left({a-1\choose2}\right)}+{\left({b+1\choose2}\right)}={\left({a\choose2}\right)}+{\left({b\choose2}\right)}-(a-1)+b<{\left({a\choose2}\right)}+{\left({b\choose2}\right)}$.
\end{pf}

\begin{prop}
We can make analogous proofs to the method described in \thmref{5pont} but for more points: in any $n$ points, there are at least $f(n)$ monochromatic edges, and if we sum the $p$'s of those edges with length different from $1$ and $d$ and it is smaller than $f(n)$, then $p_d>0$, since otherwise the expected value of the monochromatic edges would be smaller than $f(n)$. Such a method can be applied for finitely many $(1,d)$-graphs, unless the $\limsup$ of the upper bounds for the $p_x$'es for $x\neq 1$ is $0$.
\end{prop}

\begin{pfsketch}
Such a proof needs the non-edges to have a total probability less than $f(n)$, but if there are no $p_x$'es smaller than some $c>0$, then from some $n$, ${n\choose2}-2u(n)$ (that is a lower bound for the non-edges of the $(1,d)$-graph is too large and the sum of the $p_x$'s is surely larger than $f(n)$.
\end{pfsketch}

\chapter{Related problems}

\section{Spheres}

\begin{defn}
Let us suppose that for a graph $G$ there exists a point set $P(G)$ on a sphere with radius $r$ that there exists a bijection between the vertices of $G$ and the points of $P(G)$, so that the points belonging to adjacent vertices have a distance $1$ (measured in $\mathbb{R}^3$). Let us call such a graph a \nt{spherical unit distance graph} (SUDG) and also for every such radius $r$, let us call it an $r$-SUDG.
\end{defn}

\begin{defn}
Let us call such graphs \nt{faithful spherical unit distance graphs} (FSUDG) for some radius $r$, if they have realizations on a sphere with radius $r$, such that points belonging to non-adjacent vertices do not have a distance $1$ and also let us call it an $r$-FSUDG.
\end{defn}

\begin{rem}
If we would measure the distance on the sphere, we would get the same family of graphs if we take the sphere with radius ${r\over2\arcsin{1\over2r}}$ as then for the image of $P(G)$ scaled with ${r\over2\arcsin{1\over2r}}$, the spherical distances would be $1$ for all edges. The other way also works: if there would be a graph which is SUDG for the definition using spherical distances for $r\ge{1\over2}$, we could scale it back to get a SUDG with the original definition as the function ${r\over2\arcsin{1\over2r}}$ is . The same applies for FSUDG's. Therefore, we can use the original definition from here.
\end{rem}

\begin{thm}
For all $r$'s, an $r$-SUDG with $n$ vertices has at most $cn^{4\over3}$ edges.
\end{thm}

\begin{pf}
The same proof applies as in \thmref{unfelsobecsles} as all the properties we use about planar graphs can be also used on a sphere (except for that two unit circles can have the same center if $r={1\over\sqrt{2}}$, but that only changes the bound with a constant factor).
\end{pf}

\begin{defn}
Let us call a graph which has a UDR on all spheres with a large enough radius a \nt{generalized unit distance graph} (GUDG).
\end{defn}

\begin{defn}
Let us call a graph which has a FUDR on all spheres with a large enough radius a \nt{faithful generalized unit distance graph} (FGUDG).
\end{defn}

\begin{rem}
A generalized unit distance graph is not necessarily a unit distance graph.
\end{rem}

\begin{pf}
The graph in Figure \ref{gudgnemudg} is an example for a GUDG, which is not a UDG.
\end{pf}

\begin{prop}
The Moser spindle is a generalized unit distance graph.
\end{prop}

\begin{pf}
On a sufficiently large ($r\ge{1\over\sqrt{3}}$) sphere, we can take an equilateral triangle with side length $1$, since it has a circumcircle with radius $\ge{1\over\sqrt{3}}$ and we can find such a circle on the sphere. And if $r>{1\over\sqrt{3}}$, we also can find another such triangle which has a common side with the previous one, since there is a circle on the surface of the sphere, which has the same radius and intersects the previous one in two vertices of the triangle. But this way, if we take the two farthest vertices of the triangles, and take their distance $d\le\sqrt{3}$, we can draw a circle on the sphere with one of them as the center and $d$ as the radius (such a circle exists if $r\ge{\sqrt{3}\over\sqrt{2}}$). And if we intersect it with a unit circle around the other vertex, we can construct a Moser spindle, and no two points coincide for a large enough radius as their distance converges to their distance in the planar construction, whichi is positive.
\end{pf}

\begin{center}
	\input{Figures/GUDGaminemUDG}
	\captionof{figure}{}
	\label{gudgnemudg}
\end{center}

\begin{con}[description={Moser, 1966}]
{\rm{\cite{ehp}}} For some $c$, a GUDG with $n$ vertices cannot have more than $cn$ edges.
\end{con}

This conjecture turned out to be false:

\begin{thm}[description={Erdős, Hickerson, Pach, 1989}]
{\rm{\cite{ehp}}}There exist $c_1,c_2>0$ such that for every natural number $n$ and for every $0<\alpha<2$ one can find $n$ points in $S^2$ with the property that each is at a distance $\alpha$ from at least $c_1\log^{*}{n}$ others, which means that after a scaling, there is a unit distance graph on every spheres with $n$ vertices all having degree $c_1\log^{*}{n}$, so in total, it has $c\cdot n\cdot\log^{*}{n}$ edges for some $c$.
\end{thm}

\begin{pf}
Let us take a three dimensional coordinate system and a unit sphere $\lbrace (x,y,z):x^2+y^2+z^2=1\rbrace$ and let $S_{\varepsilon}$ be the set of points on the surface of this sphere, where $\abs{z}\le\varepsilon$ for any $\varepsilon\ge0$. Let us call $S_0$ the equator of $S^2$ and $S_{\varepsilon}$ as a strip of radius $\varepsilon$ around the equator. Let us fix an $\alpha$, where $0<\alpha<2$ and let us fix an $\varepsilon$, where $2\sqrt{1-\varepsilon^2}>\alpha$ to ensure that the diameter of the two circles bounding $S_{\varepsilon}$ is larger than $\alpha$.

Now we will construct a point set $P$ for all $k\ge 1$ numbers so that all points in $P$ have degree at least $k$ in the graph of $\alpha$ distances. Let us first suppose that we have such a set $P=\lbrace p_1,p_2,...,p_{n(k)}\rbrace$ for some $k$, which are inside $S_{\varepsilon(k)}$ for some $\varepsilon(k)<\varepsilon$. We also take antipodal points $u$ and $v$, where $u$ is close enough to $(0,0,1)$ (the distance will be specified later). Let us now define a rotation $\pi_i$ around $uv$ for all points $p_i$, which translates $p_i$ to a position that has distance $\alpha$ from it. It can be done because we chose $\varepsilon$ in a way that it is possible if $\delta$ is small enough. Let us now define sets $P^{(i)}$, where $P^{(0)}=P$ and $P^{(i)}=\pi_i\left(\cup_{j=0}^{i-1}{P^{(j)}}\right)$. Let us now define $P^{*}=P^{(0)}\cup_{i=1}^{n(k)}{P(i)}$. And if $\delta$ is sufficiently small, then the $\pi_i$'s exist (that is why we chose the points to be inside $S_{\varepsilon(k)}$), the $P^{(i)}$'s are disjoint (there are only finitely many possibilities for coincidence) and we still can define an $\varepsilon(k+1)$ so that all the points of $P^{*}$ are in $S_{\varepsilon(k+1)}$. And the degree of all the vertices is at least $\log^{*}{n}$.
\end{pf}

\begin{thm}[description={Swanpoel, Valtr, 2004}]
{\rm{\cite{sv}}} There exists $c>0$ such that for any $D>1$ and $n\ge 2$, $u_D(n)>cn\sqrt{\log{n}}$, where $U_D(n)$ denotes the maximal edge count of a unit distance graph on a sphere with diameter $D$.
\end{thm}

\section{Other dimensions}

The question could be asked for the graph of unit distances in $\mathbb{R}^n$ for some $n$ other than $2$. For $n=0$, the chromatic number is obviously $1$ as there is only $1$ point which can be coloured with $1$ colour and there are no edges, so this is a good colouring.

\begin{prop}
The chromatic number of $\mathbb{R}^1$ is $2$.
\end{prop}

\begin{pf}
The colouring by the parity of the integer part is good as for any $x$ coloured with red, both of its neighbours, $x-1$ and $x+1$ are coloured with blue and for any $x$ coloured with blue, both of its neighbours, $x-1$ and $x+1$ are coloured with red. And $2$ colours are needed, as for any edge, the endpoints must have different colours.
\end{pf}

\begin{thm}[description={Nechushtan, 2002}]
{\rm{\cite{n}}} The chromatic number of $\mathbb{R}^3$ is at least $6$.
\end{thm}

\begin{thm}[description={Radoičić, Tóth, }]
{\rm{\cite{rt}}} The chromatic number of $\mathbb{R}^3$ is at most $15$.
\end{thm}

\begin{defn}
Similarly to the two dimension case, we can also define $p_d$'s in three dimensions. The difference is that now we only can do it for measurable colourings.
\end{defn}

\begin{prop}
For $d\ge {1\over\sqrt{3}}$, $p_d\le{1\over3}$.
\end{prop}

\begin{pf}
For such $d$'s we can draw a tetrahedron with one side being an equilateral triangle with side length $1$, while the remaining edges having length $d$. In this graph, at most one of the $d$-edges is monochromatic, otherwise, there would be a monochromatic $1$-edge.
\end{pf}

\begin{prop}[label=2heted]
For large enough $d$, $p_d\le{2\over7}$.
\end{prop}

\begin{pf}
Let us take a sphere with radius $d$ and center $O$, and let us draw a unit Moser spindle on it. If we connect all the vertices of the Moser spindle with $O$, at most two of the $d$-edges of the graph described above are monochromatic, since at most $2$ points from those on the surface on the sphere have the same colour as $O$ because of the following lemma:

\begin{lem}
In any proper colouring of the Moser spindle, a colour occurs at most $2$ times.
\end{lem}

\begin{pf}
If we take the labelling in Figure \ref{moser}, if $A$ has the given colour, then none of $B$, $C$, $D$ and $E$ and at most one of $F$ and $G$ can have that colour, so we are done. If $A$ does not have it, then at most one of the vertices of the triangle $BCF$ and at most one of the vertices of the triangle $DEG$ can have that colour, so we are done again.
\end{pf}

And from here, the statement of Proposition \thmref{2heted} follows.
\end{pf}

\begin{prop}
For infinietly many $d$'s, there exists a $(1,d)$-graph in $\mathbb{R}^3$, which has $7$ vertices and $17$ edges.
\end{prop}

\begin{pfsketch}
First we prove that we can draw Graph 1 from \thmref{nehany1dgraf} on the surface of a sphere with center $O$ and radius $d$. We can draw two equilateral triangles with side length $1$ and exactly one common vertex for $d>{1\over\sqrt{3}}$ (since for $d\ge{1\over\sqrt{3}}$ we can draw an equilateral triangle $ABC$ with side length $1$ and if the inequality is strict, then the circumcircle of $ABC$ is not a great circle of the sphere, so there are other cirlces with the same radius through $A$ and (apart from two of them), they do not intersect circle $ABC$ in $B$ or $C$, so if we draw an equilateral triangle inside such a circle with $A$ as one of its vertices, it will satisfy the condition). And for an appropriate, but general $d$, it is possible to rotate the two triangles so that $\abs{BD}=\abs{CE}=d$. And since $BCDE$ are on a circle, they are in a plane, so we can reflect $A$ through this plane to get another point $F$, which is also connected with $B$, $C$, $D$ and $E$. So among the $7$ points, there are $17$ distances, which are $1$ or $d$: $\abs{OA}=\abs{OB}=\abs{OC}=\abs{OD}=\abs{OE}=\abs{BD}=\abs{CE}=d$ and $\abs{AB}=\abs{BC}=\abs{CA}=\abs{AD}=\abs{DE}=\abs{EA}=\abs{FB}=\abs{BC}=\abs{CF}=\abs{FD}=\abs{DE}=\abs{EF}=1$.
\end{pfsketch}

\begin{prob}
It would be interesting to find a $(1,d)$-graph in $\mathbb R^3$ with a general $d$, because (if it is for appropriate $d$'s) we could make a similar argument to Proposition \thmref{5pont}: for the lengths of the remaining segments ($x$, $y$ and $z$) we have $p_x, p_y, p_z\le{2\over7}$, so $p_d$ would be positive for an infinite number of $d$'s.
\end{prob}

\phantomsection


\addcontentsline{toc}{chapter}{Bibliography}

\begin{thebibliography}{BKT}	

\bibitem[A]{a} Ackerman, E.: On topological graphs with at most four crossings per edge https://arxiv.org/abs/1509.01932

\bibitem[BMP]{bmp}{\sc Brass, P., Moser, W. O. J., Pach, J.}: Research Problems in Discrete Geometry (2005)

\bibitem[BP]{bp}{\sc Brass, P., Pach, J.}: The Maximum Number of Times the Same Distance Can Occur among the Vertices of a Convex $n$-gon Is $O(n\log{n})$, {\it Journal of Combinatorial Theory, Series A} {\bf 94} (2001), 178--179.

\bibitem[DG]{dg}{\sc de Grey, A. D. N. J.}: The chromatic number of the plane is at least 5, {\it Geombinatorics} (2018), 28: 18-31.

\bibitem[EHP]{ehp}{\sc Erdős, P., Hickerson, D., Pach, J.}: A Problem of Leo Moser About Repeated Distances on the Sphere, {\it The American Mathematical Monthly} (1989), 96:7, 569-575.

\bibitem[F]{f}{\sc Füredi, Z.}: The maximum Number of Unit Distances in a Convex $n$-gon, {\it Journal of Combinatorial Theory, Series A} {\bf 55} (1990), 316--320.

\bibitem[GP]{gp}{\sc Globus, A., Parshall, A.}: Small unit distance graphs in the plane, https://arxiv.org/abs/1905.07829

\bibitem[H]{h}{\sc Hadwiger, H.}: Ungelöste probleme, nr. 11, {\it Elemente der Mathematik} (1961), 16: 103-104

\bibitem[He]{he}{\sc Heule, M.}: https://dustingmixon.wordpress.com/2019/07/08/polymath16-thirteenth-thread-bumping-the-deadline/\#comment-23934

\bibitem[N]{n}{\sc Nechushtan, O.}: On the space chromatic number, {\it Discrete Mathematics} (2002), 256: 499-507

\bibitem[PB1]{pb1} https://dustingmixon.wordpress.com/2018/09/14/polymath16-eleventh-thread-chromatic-numbers-of-planar-sets/

\bibitem[PB2]{pb2} https://dustingmixon.wordpress.com/2019/03/23/polymath16-twelfth-thread-year-in-review-and-future-plans/

\bibitem[PRTT]{prtt}{\sc Pach, J., Radoičić, R., Tardos, G., Tóth, G.}: Improving the crossing lemma by finding more crossings in sparse graphs, {\it Discrete and Computational Geometry} (2006), {\bf 36} (4): 527-552.

\bibitem[PW]{pw} http://michaelnielsen.org/polymath1/index.php?title=Hadwiger-Nelson\_problem

\bibitem[PWP]{pwp} http://michaelnielsen.org/polymath1/\\index.php?title=Probabilistic\_formulation\_of\_Hadwiger-Nelson\_problem

\bibitem[RT]{rt} {\sc Radoičić, R., Tóth, G.} Note on the Chromatic Number of the Space, {\it Discrete and Computational Geometry. Algorithms and Combinatorics} {\bf 25} 695--698.

\bibitem[Sch]{sch}{\sc Schade, C.}: Exakte maximale Anzahlen gleicher Abstände (1993)

\bibitem[Schae]{schae}{\sc Schaefer, M.}: Realizability of Graphs and Linkages, {\it Thirty Essays on Geometric Graph Theory} (2013), 461-482.

\bibitem[SSzT]{sszt}{\sc Spencer, J., Szemerédi, E., Trotter, W. Jr.}: Unit Distances in the Euclidean Plane, {\it Graph Theory and Combinatorics} (1984), 293--303.

\bibitem[SV]{sv}{\sc Swanpoel, K., Valtr, P.}: The Unit Distance Problem on Spheres, {\it Towards a Theory of Geometric Graphs} volume {\bf 342} of Contemporary Mathematics (2004)

\bibitem[Sz1]{sz1}{\sc Székely, L.}: Crossing numbers and hard Erdős problems in discrete geometry, {\it Combinatorics, Probability and Computing} {\bf 6} (1997), 353--358.

\bibitem[Sz2]{sz2} https://people.math.gatech.edu/\~{}hsmith90/4022/Unit\_Distances.pdf

\bibitem[T]{t} {\sc Townsend, S. P.}: Colouring the plane with no monochrome unit, {\it Geombinatorics} {\bf XIV(4)} (2005), 181-193.

\bibitem[TS]{ts} https://en.wikipedia.org/wiki/Tarski\%E2\%80\%93Seidenberg\_theorem

\bibitem[WA1]{wa1} https://www.wolframalpha.com/input/?i=maximize+(x\%2B(1-x)\%2F4)\%5E(1\%2F3)\%2B((1-x)\%2F4)\%5E(1\%2F3),+0\%3C\%3Dx\%3C\%3D1

\bibitem[WA2]{wa2} https://www.wolframalpha.com/input/?i=sqrt((d\%5E2\%2F2-1\%2F2)\%5E2\%2B(sqrt(d\%5E2-d\%5E4\%2F4)-sqrt(d\%5E2-1\%2F4))\%5E2)\%3D1\%2F2

\bibitem[WA3]{wa3} https://www.wolframalpha.com/input/?i=sqrt((d\%5E2\%2F2-1\%2F2)\%5E2\%2B(sqrt(d\%5E2-d\%5E4\%2F4)-sqrt(d\%5E2-1\%2F4))\%5E2)\%3D2\%2F(sqrt(15))

\bibitem[WA4]{wa4} https://www.wolframalpha.com/input/?i=sqrt((d\%5E2\%2F2-1\%2F2)\%5E2\%2B(sqrt(d\%5E2-d\%5E4\%2F4)-sqrt(d\%5E2-1\%2F4))\%5E2)\%3D(sqrt(3)-1)\%2F(sqrt(2))

\end{thebibliography}
\end{document}